\newcommand{\xba}{\alpha}
\newcommand{\xbb}{\beta}

\newcommand{\xbe}{\in}
\newcommand{\xbf}{\phi}
\newcommand{\xbg}{\gamma}

\newcommand{\xbm}{\mu}

\newcommand{\xbo}{\omega}

\newcommand{\xbq}{\psi}

\newcommand{\xCN}{\neg}
\newcommand{\xCQ}{\emptyset}

\newcommand{\xCd}{\approx}

\newcommand{\xCf}{\hspace{0.1em}}

\newcommand{\xcA}{\forall}

\newcommand{\xcE}{\exists}

\newcommand{\xcH}{\not\Rightarrow}
\newcommand{\xcI}{\not\Leftarrow}

\newcommand{\xcL}{\not\vdash}

\newcommand{\xcN}{\hspace{0.2em}\not\sim\hspace{-0.9em}\mid\hspace{0.8em}}

\newcommand{\xcT}{\bot}

\newcommand{\xcb}{\subset}
\newcommand{\xcc}{\subseteq}
\newcommand{\xcd}{\supseteq}
\newcommand{\xce}{\not\in}

\newcommand{\xcg}{\geq}
\newcommand{\xch}{\Rightarrow}
\newcommand{\xci}{\Leftarrow}
\newcommand{\xcj}{\Leftrightarrow}
\newcommand{\xck}{\leq}
\newcommand{\xcl}{\vdash}
\newcommand{\xcm}{\models}
\newcommand{\xcn}{\hspace{0.2em}\sim\hspace{-0.9em}\mid\hspace{0.58em}}

\newcommand{\xco}{\vee}
\newcommand{\xcp}{\rightarrow}

\newcommand{\xcr}{\leftrightarrow}
\newcommand{\xcs}{\cap}
\newcommand{\xcu}{\wedge}
\newcommand{\xcv}{\cup}

\newcommand{\xcz}{\Box}

\newcommand{\xDH}{\item }

\newcommand{\xdc}{{\cal C}}

\newcommand{\xdf}{{\cal F}}

\newcommand{\xdi}{{\cal I}}

\newcommand{\xdm}{{\cal M}}

\newcommand{\xdp}{{\cal P}}

\newcommand{\xds}{{\cal S}}

\newcommand{\xdy}{{\cal Y}}

\newcommand{\xEH}{ & }
\newcommand{\xEI}{\begin{itemize}}
\newcommand{\xEJ}{\end{itemize}}
\newcommand{\xEP}{ \\ }

\newcommand{\xEd}{\neq}
\newcommand{\xEh}{\begin{enumerate}}
\newcommand{\xEj}{\end{enumerate}}

\newcommand{\xeA}{\nabla}

\newcommand{\xFO}{\parallel}

\newcommand{\Xl}{\ldots}

\newcommand{\ol}{\overline}

\newcommand{\xssc}{\scriptsize}

\newcommand{\bl}{\begin{lemma} \rm}
\newcommand{\el}{\end{lemma}}
\newcommand{\br}{\begin{remark} \rm}
\newcommand{\er}{\end{remark}}
\newcommand{\be}{\begin{example} \rm}
\newcommand{\ee}{\end{example}}
\newcommand{\bco}{\begin{corollary} \rm}
\newcommand{\eco}{\end{corollary}}
\newcommand{\bc}{\begin{claim} \rm}
\newcommand{\ec}{\end{claim}}
\newcommand{\bfa}{\begin{fact} \rm}
\newcommand{\efa}{\end{fact}}
\newcommand{\bp}{\begin{proposition} \rm}
\newcommand{\ep}{\end{proposition}}
\newcommand{\bd}{\begin{definition} \rm}
\newcommand{\ed}{\end{definition}}
\newcommand{\bcs}{\begin{construction} \rm}
\newcommand{\ecs}{\end{construction}}
\newcommand{\bcd}{\begin{condition} \rm}
\newcommand{\ecd}{\end{condition}}
\newcommand{\bt}{\begin{theorem} \rm}
\newcommand{\et}{\end{theorem}}
\newcommand{\bn}{\begin{notation} \rm}
\newcommand{\en}{\end{notation}}
\newcommand{\bfi}{\begin{bild} \rm}
\newcommand{\efi}{\end{bild}}
\newcommand{\bsta}{\begin{statement} \rm}
\newcommand{\esta}{\end{statement}}
\newcommand{\bcom}{\begin{comment} \rm}
\newcommand{\ecom}{\end{comment}}
\newcommand{\bdia}{\begin{diagram} \rm}
\newcommand{\edia}{\end{diagram}}

\newcommand{\bfc}{\begin{figure}[htb] \begin{center}}
\newcommand{\efc}{\end{center} \end{figure}}

\documentclass{article}

\usepackage{amssymb,latexsym,epic,eepic,rotating}

\title{
Size and logic
\thanks{
Paper 339
}
}

\author{Dov M Gabbay
\thanks{
Dov.Gabbay@kcl.ac.uk, www.dcs.kcl.ac.uk/staff/dg
} \\
King's College, London
\thanks{
Department of Computer Science, King's College London, Strand,
London WC2R 2LS, UK
} \\ \\
Karl Schlechta
\thanks{
ks@cmi.univ-mrs.fr, karl.schlechta@web.de, http://www.cmi.univ-mrs.fr/ $\sim$ ks
} \\
Laboratoire d'Informatique Fondamentale de Marseille
\thanks{
UMR 6166, CNRS and Universit\'{e} de Provence,
Address: CMI, 39, rue Joliot-Curie, F-13453 Marseille Cedex 13, France
}
}

\begin{document}

\newtheorem{lemma}{Lemma}[section]
\newtheorem{theorem}[lemma]{Theorem}
\newtheorem{proposition}[lemma]{Proposition}
\newtheorem{corollary}[lemma]{Corollary}
\newtheorem{claim}[lemma]{Claim}
\newtheorem{fact}[lemma]{Fact}
\newtheorem{remark}[lemma]{Remark}
\newtheorem{definition}{Definition}[section]
\newtheorem{construction}{Construction}[section]
\newtheorem{condition}{Condition}[section]
\newtheorem{example}{Example}[section]
\newtheorem{notation}{Notation}[section]
\newtheorem{bild}{Figure}[section]
\newtheorem{comment}{Comment}[section]
\newtheorem{statement}{Statement}[section]
\newtheorem{diagram}{Diagram}[section]

\renewcommand{\labelenumi}
  {(\arabic{enumi})}
\renewcommand{\labelenumii}
  {(\arabic{enumi}.\arabic{enumii})}
\renewcommand{\labelenumiii}
  {(\arabic{enumi}.\arabic{enumii}.\arabic{enumiii})}
\renewcommand{\labelenumiv}
  {(\arabic{enumi}.\arabic{enumii}.\arabic{enumiii}.\arabic{enumiv})}

\maketitle

\begin{abstract}

We show how to develop a multitude of rules of nonmonotonic logic
from very simple and natural notions of size, using them as
building blocks.

\end{abstract}

\tableofcontents

\setcounter{secnumdepth}{3}
\setcounter{tocdepth}{3}

%
%
%
\section{
Introduction
}
\subsection{
Context
}

The study of modal and temporal logic and the study of
substructural logic went for many years along the following lines:
on the one hand we had syntactic proof theoretic systems, mainly Gentzen
or Hilbert systems and on the other hand we had semantical
interpretations, mainly possible worlds or algebraic structures, and
the community very thoroughly analysed properties of one against
matching properties of the other.
The success of such depended on the correct identification of the
correct features in the semantics.
In the case of nonmonotonic logic there is the syntactical consequence
relation on the one hand and the preferential ordering on the other but
the research is not yet in a similar comprehensive stage as in the other
areas.
In this paper we use the important semantical feature of size
to display a detailed matching between syntactical and
semantical conditions for non monotonic systems.

We show how one can develop a multitude of rules for nonmonotonic
logics from a very small
set of principles about reasoning with size. The notion of size gives an
algebraic semantics to nonmonotonic logics, in the sense that $ \xba $
implies $ \xbb $
iff the set of cases where $ \xba \xcu \xCN \xbb $ holds is a small subset
of all $ \xba -$cases.
In a similar way, e.g. Heyting algebras are an algebraic semantics for
intuitionistic logic.

In our understanding, algebraic semantics describe
the abstract properties corresponding model sets have. Structural
semantics, on the other hand, give intuitive concepts like accessibility
or preference, from which properties of model sets, and thus algebraic
semantics, originate.

Varying properties of structural semantics (e.g. transitivity, etc.)
result in
varying properties of algebraic semantics, and thus of logical rules.
We consider operations directly on the algebraic semantics and their
logical
consequences, and we see that simple manipulations of the size concept
result in most rules of nonmonotonic logics. Even more, we show how
to generate new rules from those manipulations. The result is one big
table, which, in a much more modest scale, can be seen as a
``periodic table'' of the ``elements'' of nonmonotonic logic.
Some simple underlying principles allow to generate them all.

Historical remarks: The first time that abstract size was related to
nonmonotonic logics was, to our knowledge,
in the second author's  \cite{Sch90} and  \cite{Sch95-1}, and,
independently, in  \cite{BB94}.
More detailed remarks can e.g. be found in  \cite{GS08c}. But,
again to our knowledge, connections are elaborated systematically and in
fine
detail here for the first time.
\subsection{
Overview
}

The main part of this paper is the big table
in Section \ref{Section Main-Table} (page \pageref{Section Main-Table}).
It shows connections and how to develop a multitude of logical rules
known from nonmonotonic logics by combining a small number of principles
about size. We use them as building blocks to construct the rules from.

These principles are some basic and very natural postulates,
$ \xCf (Opt),$ $ \xCf (iM),$ $(eM \xdi ),$ $(eM \xdf ),$ and a continuum
of power of the notion of
``small'', or, dually, ``big'', from $(1*s)$ to $(< \xbo *s).$
From these, we can develop the rest except, essentially, Rational
Monotony,
and thus an infinity of different rules.

This is a conceptual paper, and it does not contain any more difficult
formal results. The interest lies, in our opinion, in the simplicity,
paucity,
and naturalness of the basic building blocks.
We hope that this schema brings more and deeper order into the rich fauna
of nonmonotonic and related logics.
\section{
Main table
}

$\hspace{0.01em}$

{\xssc LABEL: {Section Table}}
\label{Section Table}
\subsection{
Notation
}

 \xEh

 \xDH

$ \xdp (X)$ is the power set of $X,$ $ \xcc $ is the subset relation, $
\xcb $ the strict part of
$ \xcc,$ i.e. $A \xcb B$ iff $A \xcc B$ and $A \xEd B.$
The operators $ \xcu,$ $ \xCN,$ $ \xco,$ $ \xcp $ and $ \xcl $ have
their usual, classical interpretation.

 \xDH

$ \xdi (X) \xcc \xdp (X)$ and $ \xdf (X) \xcc \xdp (X)$ are dual abstract
notions of size, $ \xdi (X)$ is the
set of ``small'' subsets of $X,$ $ \xdf (X)$ the set of ``big'' subsets of
$X.$ They are
dual in the sense that $A \xbe \xdi (X) \xcj X-A \xbe \xdf (X).$ `` $ \xdi
$ '' evokes ``ideal'',
`` $ \xdf $ '' evokes ``filter'' though the full strength of both is reached
only
in $(< \xbo *s).$ ``s'' evokes ``small'', and `` $(x*s)$ '' stands for
`` $x$ small sets together are still not everything''.

 \xDH

If $A \xcc X$ is neither in $ \xdi (X),$ nor in $ \xdf (X),$ we say it has
medium size, and
we define $ \xdm (X):= \xdp (X)-( \xdi (X) \xcv \xdf (X)).$ $
\xdm^{+}(X):= \xdp (X)- \xdi (X)$ is the set of subsets
which are not small.

 \xDH

$ \xeA x \xbf $ is a generalized first order quantifier, it is read
``almost all $x$ have property $ \xbf $ ''. $ \xeA x( \xbf: \xbq )$ is the
relativized version, read:
``almost all $x$ with property $ \xbf $ have also property $ \xbq $ ''. To
keep the table
simple, we write mostly only the non-relativized versions.
Formally, we have $ \xeA x \xbf: \xcj \{x: \xbf (x)\} \xbe \xdf (U)$
where $U$ is the universe, and
$ \xeA x( \xbf: \xbq ): \xcj \{x:( \xbf \xcu \xbq )(x)\} \xbe \xdf (\{x:
\xbf (x)\}).$
Soundness and completeness results on $ \xeA $ can be found in
 \cite{Sch95-1}.

 \xDH

Analogously, for propositional logic, we define:

$ \xba \xcn \xbb $ $: \xcj $ $M( \xba \xcu \xbb ) \xbe \xdf (M( \xba )),$

where $M( \xbf )$ is the set of models of $ \xbf.$

 \xDH

In preferential structures, $ \xbm (X) \xcc X$ is the set of minimal
elements of $X.$
This generates a principal filter by $ \xdf (X):=\{A \xcc X: \xbm (X) \xcc
A\}.$ Corresponding
properties about $ \xbm $ are not listed systematically.

 \xDH

The usual rules $ \xCf (AND)$ etc. are named here $(AND_{ \xbo }),$ as
they are in a
natural ascending line of similar rules, based on strengthening of the
filter/ideal properties.

 \xEj
\subsection{
The groupes of rules
}

The rules are divided into 5 groups:

 \xEh

 \xDH $ \xCf (Opt),$ which says that ``All'' is optimal - i.e. when there
are no
exceptions, then a soft rule $ \xcn $ holds.

 \xDH 3 monotony rules:

 \xEh
 \xDH $ \xCf (iM)$ is inner monotony, a subset of a small set is small,
 \xDH $(eM \xdi )$ external monotony for ideals: enlarging the base set
keeps small
sets small,
 \xDH $(eM \xdf )$ external monotony for filters: a big subset stays big
when the base
set shrinks.
 \xEj

These three rules are very natural if ``size'' is anything coherent over
change
of base sets. In particular, they can be seen as weakening.

 \xDH $( \xCd )$ keeps proportions, it is here mainly to point the
possibility out.

 \xDH a group of rules $x*s,$ which say how many small sets will not yet
add to
the base set.
 \xDH Rational monotony, which can best be understood as robustness of $
\xdm^{+},$
see $( \xdm^{++})(3).$

 \xEj
\subsubsection{
Regularities
}

 \xEh

 \xDH

The group of rules $(x*s)$ use ascending strength of $ \xdi / \xdf.$

 \xDH

The column $( \xdm^{+})$ contains interesting algebraic properties. In
particular,
they show a strengthening from $(3*s)$ up to Rationality. They are not
necessarily
equivalent to the corresponding $(I_{x})$ rules, not even in the presence
of the basic rules. The examples show that care has to be taken when
considering
the different variants.

 \xDH

Adding the somewhat superflous $(CM_{2}),$ we have increasing cautious
monotony from $ \xCf (wCM)$ to full $(CM_{ \xbo }).$

 \xDH

We have increasing ``or'' from $ \xCf (wOR)$ to full $(OR_{ \xbo }).$

 \xDH

The line $(2*s)$ is only there because there seems to be no $(
\xdm^{+}_{2}),$ otherwise
we could begin $(n*s)$ at $n=2.$

 \xEj
\subsection{
Direct correspondences
}

Several correspondences are trivial and are mentioned now.
Somewhat less obvious (in)dependencies are given in
Section \ref{Section Coherent-Systems} (page \pageref{Section Coherent-Systems})
.
Finally, the connections with the $ \xbm -$rules are given in
Section \ref{Section Principal} (page \pageref{Section Principal}).
In those rules, $(I_{ \xbo })$ is implicit, as they are about principal
filters.
Still, the $ \xbm -$rules are written in the main table in their
intuitively
adequate place.

 \xEh

 \xDH

The columns ``Ideal'' and ``Filter'' are mutually dual, when both entries are
defined.

 \xDH

The correspondence between the ideal/filter column and the $ \xeA -$column
is obvious,
the latter is added only for completeness' sake, and to point out the
trivial translation to first order logic.

 \xDH

The ideal/filter and the AND-column correspond directly.

 \xDH

We can construct logical rules from the $ \xdm^{+}-column$ by direct
correspondence,
e.g. for $( \xdm^{+}_{ \xbo }),$ (1):

Set $Y:=M( \xbg ),$ $X:=M( \xbg \xcu \xbb ),$ $A:=M( \xbg \xcu \xbb \xcu
\xba ).$

 \xEI
 \xDH $X \xbe \xdm^{+}(Y)$ will become $ \xbg \xcN \xCN \xbb $
 \xDH $A \xbe \xdf (X)$ will become $ \xbg \xcu \xbb \xcn \xba $
 \xDH $A \xbe \xdm^{+}(Y)$ will become $ \xbg \xcN \xCN ( \xba \xcu \xbb
).$
 \xEJ

so we obtain $ \xbg \xcN \xCN \xbb,$ $ \xbg \xcu \xbb \xcn \xba $ $ \xch
$ $ \xbg \xcN \xCN ( \xba \xcu \xbb ).$

We did not want to make the table too complicated, so
such rules are not listed in the table.

 \xDH

Various direct correspondences:

 \xEI

 \xDH In the line $ \xCf (Opt),$ the filter/ideal entry corresponds to $
\xCf (SC),$
 \xDH in the line $ \xCf (iM),$ the filter/ideal entry corresponds to $
\xCf (RW),$
 \xDH in the line $(eM \xdi ),$ the ideal entry corresponds to $(PR' )$
and $ \xCf (wOR),$
 \xDH in the line $(eM \xdf ),$ the filter entry corresponds to $ \xCf
(wCM),$
 \xDH in the line $( \xCd ),$ the filter/ideal entry corresponds to $ \xCf
(disjOR),$
 \xDH in the line $(1*s),$ the filter/ideal entry corresponds to $ \xCf
(CP),$
 \xDH in the line $(2*s),$ the filter/ideal entry corresponds to
$(CM_{2})=(OR_{2}).$

 \xEJ

 \xDH

Note that one can, e.g., write $(AND_{2})$ in two flavours:

 \xEI
 \xDH $ \xba \xcn \xbb,$ $ \xba \xcn \xbb ' $ $ \xch $ $ \xba \xcL \xCN
\xbb \xco \xCN \xbb ',$ or
 \xDH $ \xba \xcn \xbb $ $ \xch $ $ \xba \xcN \xCN \xbb $
 \xEJ

(which is $(CM_{2})=(OR_{2}).)$

For reasons of simplicity, we mention only one.

 \xEj
\subsection{
Rational Monotony
}

$ \xCf (RatM)$ does not fit into adding small sets. We have exhausted the
combination
of small sets by $(< \xbo *s),$ unless we go to languages with infinitary
formulas.

The next idea would be to add medium size sets. But, by definition,
$2*medium$
can be all. Adding small and medium sets would not help either: Suppose we
have a rule $medium+n*small \xEd all.$ Taking the complement of the first
medium
set, which is again medium, we have the rule $2*n*small \xEd all.$ So we
do not
see any meaningful new internal rule. i.e. without changing the base set.

Probably, $ \xCf (RatM)$ has more to do with independence: by default, all
``normalities'' are independent, and intersecting with another formula
preserves normality.
\subsection{
Summary
}

We can obtain all rules except $ \xCf (RatM)$ and $( \xCd )$ from $ \xCf
(Opt),$ the monotony
rules - $ \xCf (iM),$ $(eM \xdi ),$ $(eM \xdf )$ -, and $(x*s)$ with
increasing $x.$
\subsection{
Main table
}

$\hspace{0.01em}$

{\xssc LABEL: {Section Main-Table}}
\label{Section Main-Table}
\newpage

\begin{turn}{90}

{\tiny

\begin{tabular}{|c|c@{:}c|c|c|c|c|c|c|}

\hline

\xEH
``Ideal''
\xEH
``Filter''
\xEH
$ \xdm^+ $
\xEH
$ \xeA $
\xEH
various rules
\xEH
AND
\xEH
OR
\xEH
Caut./Rat.Mon.
\xEP

\hline
\hline

\multicolumn{9}{|c|}{Optimal proportion} \xEP

\hline

$(Opt)$
\xEH
$ \xCQ \xbe \xdi (X)$
\xEH
$X \xbe \xdf (X)$
\xEH
\xEH
$ \xcA x \xbf \xcp \xeA x \xbf$
\xEH
$(SC)$
\xEH
\xEH
\xEH
\xEP

\xEH
\xEH
\xEH
\xEH
\xEH
$ \xba \xcl \xbb \xch \xba \xcn \xbb $
\xEH
\xEH
\xEH
\xEP

\hline
\hline

\multicolumn{9}{|c|}
{Monotony (Improving proportions). $(iM)$: internal monotony,
$(eM \xdi )$: external monotony for ideals,
$(eM \xdf )$: external monotony for filters}
\xEP

\hline

$(iM)$
\xEH
$A \xcc B \xbe \xdi (X)$ $ \xch $
\xEH
$A \xbe \xdf (X)$, $A \xcc B \xcc X$
\xEH
\xEH
$\xeA x \xbf \xcu \xcA x (\xbf \xcp \xbf')$
\xEH
$(RW)$
\xEH
\xEH
\xEH
\xEP

\xEH
$A \xbe \xdi (X)$
\xEH
$ \xch $ $B \xbe \xdf (X)$
\xEH
\xEH
$ \xcp $ $ \xeA x \xbf'$
\xEH
$ \xba \xcn \xbb, \xbb \xcl \xbb' \xch $
\xEH
\xEH
\xEH
\xEP

\xEH
\xEH
\xEH
\xEH
\xEH
$ \xba \xcn \xbb' $
\xEH
\xEH
\xEH
\xEP

\hline

$(eM \xdi )$
\xEH
$X \xcc Y \xch$
\xEH
\xEH
\xEH
$\xeA x (\xbf: \xbq) \xcu$
\xEH
$(PR')$
\xEH
\xEH
$(wOR)$
\xEH
\xEP

\xEH
$\xdi (X) \xcc \xdi (Y)$
\xEH
\xEH
\xEH
$\xcA x (\xbf' \xcp \xbq) \xcp$
\xEH
$\xba \xcn \xbb, \xba \xcl \xba',$
\xEH
\xEH
$ \xba \xcn \xbb,$ $ \xba ' \xcl \xbb $ $ \xch $
\xEH
\xEP

\xEH
\xEH
\xEH
\xEH
$\xeA x (\xbf \xco \xbf': \xbq)$
\xEH
$\xba' \xcu \xCN \xba \xcl \xbb \xch$
\xEH
\xEH
$ \xba \xco \xba ' \xcn \xbb $
\xEH
\xEP

\xEH
\xEH
\xEH
\xEH
\xEH
$\xba' \xcn \xbb$
\xEH
\xEH
$(\xbm wOR)$
\xEH
\xEP

\xEH
\xEH
\xEH
\xEH
\xEH
$(\xbm PR)$
\xEH
\xEH
$\xbm(X \xcv Y) \xcc \xbm(X) \xcv Y$
\xEH
\xEP

\xEH
\xEH
\xEH
\xEH
\xEH
$X \xcc Y \xch$
\xEH
\xEH
\xEH
\xEP

\xEH
\xEH
\xEH
\xEH
\xEH
$\xbm(Y) \xcs X \xcc \xbm(X)$
\xEH
\xEH
\xEH
\xEP

\hline

$(eM \xdf )$
\xEH
\xEH
$X \xcc Y \xch$
\xEH
\xEH
$\xeA x (\xbf: \xbq) \xcu$
\xEH
\xEH
\xEH
\xEH
$(wCM)$
\xEP

\xEH
\xEH
$\xdf (Y) \xcs \xdp (X) \xcc \xdf (X)$
\xEH
\xEH
$\xcA x (\xbq \xcu \xbf \xcp \xbf') \xcp$
\xEH
\xEH
\xEH
\xEH
$\xba \xcn \xbb, \xba' \xcl \xba,$
\xEP

\xEH
\xEH
\xEH
\xEH
$\xeA x (\xbf \xcu \xbf': \xbq)$
\xEH
\xEH
\xEH
\xEH
$\xba \xcu \xbb \xcl \xba' \xch$
\xEP

\xEH
\xEH
\xEH
\xEH
\xEH
\xEH
\xEH
\xEH
$\xba' \xcn \xbb$
\xEP

\hline
\hline

\multicolumn{9}{|c|}{Keeping proportions} \xEP

\hline

$(\xCd)$
\xEH
$(\xdi \xcv disj)$
\xEH
$(\xdf \xcv disj)$
\xEH
\xEH
$\xeA x(\xbf: \xbq) \xcu$
\xEH
\xEH
\xEH
$(disjOR)$
\xEH
\xEP

\xEH
$A \xbe \xdi (X),$ $B \xbe \xdi (Y),$
\xEH
$A \xbe \xdf (X),$ $B \xbe \xdf (Y),$
\xEH
\xEH
$\xeA x(\xbf': \xbq) \xcu$
\xEH
\xEH
\xEH
$ \xbf \xcn \xbq,$ $ \xbf ' \xcn \xbq ' $
\xEH
\xEP

\xEH
$X \xcs Y= \xCQ $ $ \xch $
\xEH
$X \xcs Y= \xCQ $ $ \xch $
\xEH
\xEH
$\xCN \xcE x(\xbf \xcu \xbf') \xcp$
\xEH
\xEH
\xEH
$ \xbf \xcl \xCN \xbf ',$ $ \xch $
\xEH
\xEP

\xEH
$A \xcv B \xbe \xdi (X \xcv Y)$
\xEH
$A \xcv B \xbe \xdf (X \xcv Y)$
\xEH
\xEH
$\xeA x(\xbf \xco \xbf': \xbq)$
\xEH
\xEH
\xEH
$ \xbf \xco \xbf ' \xcn \xbq \xco \xbq ' $
\xEH
\xEP

\xEH
\xEH
\xEH
\xEH
\xEH
\xEH
\xEH
$(\xbm disjOR)$
\xEH
\xEP

\xEH
\xEH
\xEH
\xEH
\xEH
\xEH
\xEH
$X \xcs Y = \xCQ \xch$
\xEH
\xEP

\xEH
\xEH
\xEH
\xEH
\xEH
\xEH
\xEH
$\xbm(X \xcv Y) \xcc \xbm(X) \xcv \xbm(Y)$
\xEH
\xEP

\hline
\hline

\multicolumn{9}{|c|}{Robustness of proportions: $n*small \xEd All$} \xEP

\hline

$(1*s)$
\xEH
$(\xdi_1)$
\xEH
$(\xdf_1)$
\xEH
\xEH
$(\xeA_1)$
\xEH
$(CP)$
\xEH
$(AND_{1})$
\xEH
\xEH
\xEP

\xEH
$X \xce \xdi (X)$
\xEH
$ \xCQ \xce \xdf (X)$
\xEH
\xEH
$ \xeA x \xbf \xcp \xcE x \xbf $
\xEH
$\xbf\xcn\xcT \xch \xbf\xcl\xcT$
\xEH
$ \xba \xcn \xbb $ $ \xch $ $ \xba \xcL \xCN \xbb $
\xEH
\xEH
\xEP

\hline

$(2*s)$
\xEH
$(\xdi_2)$
\xEH
$(\xdf_2)$
\xEH
\xEH
$(\xeA_2)$
\xEH
\xEH
$(AND_{2})$
\xEH
$(OR_{2})$
\xEH
$(CM_{2})$
\xEP

\xEH
$A,B \xbe \xdi (X) \xch $
\xEH
$A,B \xbe \xdf (X) \xch $
\xEH
\xEH
$ \xeA x \xbf \xcu \xeA x \xbq $
\xEH
\xEH
$ \xba \xcn \xbb,$ $ \xba \xcn \xbb ' $ $ \xch $
\xEH
$ \xba \xcn \xbb \xch \xba \xcN \xCN \xbb $
\xEH
$ \xba \xcn \xbb \xch \xba \xcN \xCN \xbb $
\xEP

\xEH
$ A \xcv B \xEd X$
\xEH
$A \xcs B \xEd \xCQ $
\xEH
\xEH
$ \xcp $ $ \xcE x( \xbf \xcu \xbq )$
\xEH
\xEH
$ \xba \xcL \xCN \xbb \xco \xCN \xbb ' $
\xEH
\xEH
\xEP

\hline

%
%
%

$(n*s)$
\xEH
$(\xdi_n)$
\xEH
$(\xdf_n)$
\xEH
$( \xdm^{+}_{n})$
\xEH
$(\xeA_n)$
\xEH
\xEH
$(AND_{n})$
\xEH
$(OR_{n})$
\xEH
$(CM_{n})$
\xEP

$(n \xcg 3)$
\xEH
$A_{1},.,A_{n} \xbe \xdi (X) $
\xEH
$A_{1},.,A_{n} \xbe \xdi (X) $
\xEH
$X_{1} \xbe \xdf (X_{2}),., $
\xEH
$ \xeA x \xbf_{1} \xcu.  \xcu \xeA x \xbf_{n} $
\xEH
\xEH
$ \xba \xcn \xbb_{1},., \xba \xcn \xbb_{n}$ $ \xch $
\xEH
$ \xba_{1} \xcn \xbb,., \xba_{n-1} \xcn \xbb $
\xEH
$ \xba \xcn \xbb_{1},., \xba \xcn \xbb_{n-1}$
\xEP

\xEH
$ \xch $
\xEH
$ \xch $
\xEH
$ X_{n-1} \xbe \xdf (X_{n})$ $ \xch $
\xEH
$ \xcp $
\xEH
\xEH
$ \xba \xcL \xCN \xbb_{1} \xco.  \xco \xCN \xbb_{n}$
\xEH
$ \xch $
\xEH
$ \xch $
\xEP

\xEH
$ A_{1} \xcv.  \xcv A_{n} \xEd X $
\xEH
$A_{1} \xcs.  \xcs A_{n} \xEd \xCQ$
\xEH
$X_{1} \xbe \xdm^{+}(X_{n})$
\xEH
$ \xcE x (\xbf_{1} \xcu.  \xcu \xbf_{n}) $
\xEH
\xEH
\xEH
$ \xba_{1} \xco.  \xco \xba_{n-1} \xcN \xCN \xbb $
\xEH
$ \xba \xcu \xbb_1 \xcu.  \xcu \xbb_{n-2} \xcN \xCN \xbb_{n-1}$
\xEP

\hline

$(< \xbo*s)$
\xEH
$(\xdi_\xbo)$
\xEH
$(\xdf_\xbo)$
\xEH
$( \xdm^{+}_{ \xbo })$
\xEH
$(\xeA_{\xbo})$
\xEH
\xEH
$(AND_{ \xbo })$
\xEH
$(OR_{ \xbo })$
\xEH
$(CM_{ \xbo })$
\xEP

\xEH
$A,B \xbe \xdi (X) \xch $
\xEH
$A,B \xbe \xdf (X) \xch $
\xEH
(1)
\xEH
$ \xeA x \xbf \xcu \xeA x \xbq \xcp $
\xEH
\xEH
$ \xba \xcn \xbb,$ $ \xba \xcn \xbb ' $ $ \xch $
\xEH
$ \xba \xcn \xbb,$ $ \xba ' \xcn \xbb $ $ \xch $
\xEH
$ \xba \xcn \xbb,$ $ \xba \xcn \xbb ' $ $ \xch $
\xEP

\xEH
$ A \xcv B \xbe \xdi (X)$
\xEH
$ A \xcs B \xbe \xdf (X)$
\xEH
$A \xbe \xdf (X),$ $X \xbe \xdm^{+}(Y)$
\xEH
$ \xeA x( \xbf \xcu \xbq )$
\xEH
\xEH
$ \xba \xcn \xbb \xcu \xbb ' $
\xEH
$ \xba \xco \xba ' \xcn \xbb $
\xEH
$ \xba \xcu \xbb \xcn \xbb ' $
\xEP

\xEH
\xEH
\xEH
$ \xch $ $A \xbe \xdm^{+}(Y)$
\xEH
\xEH
\xEH
\xEH
$(\xbm OR)$
\xEH
$(\xbm CM)$
\xEP

\xEH
\xEH
\xEH
(2)
\xEH
\xEH
\xEH
\xEH
$\xbm(X \xcv Y) \xcc \xbm(X) \xcv \xbm(Y)$
\xEH
$\xbm(X) \xcc Y \xcc X \xch$
\xEP

\xEH
\xEH
\xEH
$A \xbe \xdm^{+}(X),$ $X \xbe \xdf (Y)$
\xEH
\xEH
\xEH
\xEH
\xEH
$\xbm(Y) \xcc \xbm(X)$
\xEP

\xEH
\xEH
\xEH
$ \xch $ $A \xbe \xdm^{+}(Y)$
\xEH
\xEH
\xEH
\xEH
\xEH
\xEP

\xEH
\xEH
\xEH
(3)
\xEH
\xEH
\xEH
\xEH
\xEH
\xEP

\xEH
\xEH
\xEH
$A \xbe \xdf (X),$ $X \xbe \xdf (Y)$
\xEH
\xEH
\xEH
\xEH
\xEH
\xEP

\xEH
\xEH
\xEH
$ \xch $ $A \xbe \xdf (Y)$
\xEH
\xEH
\xEH
\xEH
\xEH
\xEP

\xEH
\xEH
\xEH
(4)
\xEH
\xEH
\xEH
\xEH
\xEH
\xEP

\xEH
\xEH
\xEH
$A,B \xbe \xdi (X)$ $ \xch $
\xEH
\xEH
\xEH
\xEH
\xEH
\xEP

\xEH
\xEH
\xEH
$A-B \xbe \xdi (X-$B)
\xEH
\xEH
\xEH
\xEH
\xEH
\xEP

\hline
\hline

\multicolumn{9}{|c|}{Robustness of $\xdm^+$} \xEP

\hline

$(\xdm^{++})$
\xEH
\xEH
\xEH
$(\xdm^{++})$
\xEH
\xEH
\xEH
\xEH
\xEH
$(RatM)$
\xEP

\xEH
\xEH
\xEH
(1)
\xEH
\xEH
\xEH
\xEH
\xEH
$ \xbf \xcn \xbq,  \xbf \xcN \xCN \xbq '   \xch $
\xEP

\xEH
\xEH
\xEH
$A \xbe \xdi (X),$ $B \xce \xdf (X)$
\xEH
\xEH
\xEH
\xEH
\xEH
$ \xbf \xcu \xbq ' \xcn \xbq $
\xEP

\xEH
\xEH
\xEH
$ \xch $ $A-B \xbe \xdi (X-B)$
\xEH
\xEH
\xEH
\xEH
\xEH
$(\xbm RatM)$
\xEP

\xEH
\xEH
\xEH
(2)
\xEH
\xEH
\xEH
\xEH
\xEH
$X \xcc Y,$
\xEP

\xEH
\xEH
\xEH
$A \xbe \xdf (X), B \xce \xdf (X)$
\xEH
\xEH
\xEH
\xEH
\xEH
$X \xcs \xbm(Y) \xEd \xCQ \xch$
\xEP

\xEH
\xEH
\xEH
$ \xch $ $A-B \xbe \xdf (X-$B)
\xEH
\xEH
\xEH
\xEH
\xEH
$\xbm(X) \xcc \xbm(Y) \xcs X$
\xEP

\xEH
\xEH
\xEH
(3)
\xEH
\xEH
\xEH
\xEH
\xEH
\xEP

\xEH
\xEH
\xEH
$A \xbe \xdm^+ (X),$
\xEH
\xEH
\xEH
\xEH
\xEH
\xEP

\xEH
\xEH
\xEH
$X \xbe \xdm^+ (Y)$
\xEH
\xEH
\xEH
\xEH
\xEH
\xEP

\xEH
\xEH
\xEH
$ \xch $ $A \xbe \xdm^+ (Y)$
\xEH
\xEH
\xEH
\xEH
\xEH
\xEP

\hline

\hline

\end{tabular}

}

\end{turn}

\newpage
\section{
Coherent systems
}

$\hspace{0.01em}$

{\xssc LABEL: {Section Coherent-Systems}}
\label{Section Coherent-Systems}
\subsection{
Definition and basic facts
}

Note that whenever we work with model sets, the rule

$ \xCf (LLE),$ left logical equivalence, $ \xcl \xba \xcr \xba ' $ $ \xch
$ $( \xba \xcn \xbb $ $ \xcj $ $ \xba ' \xcn \xbb )$

will hold. We will not mention this any further.

\bd

$\hspace{0.01em}$

(+++ Orig. No.:  Definition CoherentSystem +++)

$\hspace{0.01em}$

{\xssc LABEL: {Definition CoherentSystem}}
\label{Definition CoherentSystem}

A coherent system of sizes, $ \xdc \xds,$ consists of a universe $U,$ $
\xCQ \xce \xdy \xcc \xdp (U),$ and
for all $X \xbe \xdy $ a system $ \xdi (X) \xcc \xdp (X)$ (dually $ \xdf
(X),$ i.e. $A \xbe \xdf (X) \xcj X-A \xbe \xdi (X)).$
$ \xdy $ may satisfy certain closure properties like closure under $ \xcv
,$ $ \xcs,$
complementation, etc. We will mention this when needed, and not obvious.

We say that $ \xdc \xds $ satisfies a certain property iff all $X,Y \xbe
\xdy $ satisfy this
property.

$ \xdc \xds $ is called basic or level 1 iff it satisfies $ \xCf (Opt),$ $
\xCf (iM),$ $(eM \xdi ),$
$(eM \xdf ),$ $(1*s).$

$ \xdc \xds $ is level $x$ iff it satisfies $ \xCf (Opt),$ $ \xCf (iM),$
$(eM \xdi ),$ $(eM \xdf ),$ $(x*s).$

\ed

\bfa

$\hspace{0.01em}$

(+++ Orig. No.:  Fact 1-element +++)

$\hspace{0.01em}$

{\xssc LABEL: {Fact 1-element}}
\label{Fact 1-element}

Note that, if for any $Y$ $ \xdi (Y)$ consists only of subsets of at most
1 element,
then $(eM \xdf )$ is trivially satisfied for $Y$ and its subsets by $ \xCf
(Opt).$ $ \xcz $
\\[3ex]

\efa

\bfa

$\hspace{0.01em}$

(+++ Orig. No.:  Fact Not-2*s +++)

$\hspace{0.01em}$

{\xssc LABEL: {Fact Not-2*s}}
\label{Fact Not-2*s}

Let a $ \xdc \xds $ be given s.t. $ \xdy = \xdp (U).$ If $X \xbe \xdy $
satisfies $( \xdm^{++}),$ but not
$(< \xbo *s),$ then there is $Y \xbe \xdy $ which does not satisfy
$(2*s).$

\efa

\subparagraph{
Proof
}

$\hspace{0.01em}$

(+++ Orig.:  Proof +++)

We work with version (1) of $( \xdm^{++}),$ we will see in
Fact \ref{Fact M-plus-plus} (page \pageref{Fact M-plus-plus})  that all three
versions are equivalent.

As $X$ does not satisfy $(< \xbo *s),$ there are $A,B \xbe \xdi (X)$ s.t.
$A \xcv B \xbe \xdm^{+}(X).$
$A \xbe \xdi (X),$ $A \xcv B \xbe \xdm^{+}(X)$ $ \xch $ $X-(A \xcv B) \xce
\xdf (X),$ so by $( \xdm^{++})(1)$
$A=A-(X-(A \xcv B)) \xbe \xdi (X-(X-(A \xcv B)))= \xdi (A \xcv B).$
Likewise $B \xbe \xdi (A \xcv B),$
so $(2*s)$ does not hold for $A \xcv B.$ $ \xcz $
\\[3ex]

\bfa

$\hspace{0.01em}$

(+++ Orig. No.:  Fact Independence-eM +++)

$\hspace{0.01em}$

{\xssc LABEL: {Fact Independence-eM}}
\label{Fact Independence-eM}

$(eM \xdi )$ and $(eM \xdf )$ are formally independent, though intuitively
equivalent.

\efa

\subparagraph{
Proof
}

$\hspace{0.01em}$

(+++ Orig.:  Proof +++)

Let $U:=\{x,y,z\},$ $X:=\{x,z\},$ $ \xdy:= \xdp (U)-\{ \xCQ \}$

(1) Let $ \xdf (U):=\{A \xcc U:z \xbe A\},$ $ \xdf (Y)=\{Y\}$ for all $Y
\xcb U.$
$ \xCf (Opt),$ $ \xCf (iM)$ hold, $(eM \xdi )$ holds trivially, so does
$(< \xbo *s),$
but $(eM \xdf )$ fails for $U$ and $X.$

(2) Let $ \xdf (X):=\{\{z\},X\},$ $ \xdf (Y):=\{Y\}$ for all $Y \xcc U,$
$Y \xEd X.$
$ \xCf (Opt),$ $ \xCf (iM),$ $(< \xbo *s)$ hold trivially, $(eM \xdf )$
holds by
Fact \ref{Fact 1-element} (page \pageref{Fact 1-element}). $(eM \xdi )$ fails,
as $\{x\} \xbe \xdi
(X),$ but $\{x\} \xce \xdi (U).$

$ \xcz $
\\[3ex]

\bfa

$\hspace{0.01em}$

(+++ Orig. No.:  Fact Level-n-n+1 +++)

$\hspace{0.01em}$

{\xssc LABEL: {Fact Level-n-n+1}}
\label{Fact Level-n-n+1}

A level $n$ system is strictly weaker than a level $n+1$ system.

\efa

\subparagraph{
Proof
}

$\hspace{0.01em}$

(+++ Orig.:  Proof +++)

Consider $U:=\{1, \Xl,n+1\},$ $ \xdy:= \xdp (U)-\{ \xCQ \}.$ Let $ \xdi
(U):=\{ \xCQ \} \xcv \{\{x\}:x \xbe U\},$
$ \xdi (X):=\{ \xCQ \}$ for $X \xEd U.$
$ \xCf (iM),$ $(eM \xdi ),$ $(eM \xdf )$ hold trivially.
$(n*s)$ holds trivially for $X \xEd U,$ but also for $U.$ $((n+1)*s)$ does
not hold
for $U.$ $ \xcz $
\\[3ex]

\br

$\hspace{0.01em}$

(+++ Orig. No.:  Remark Infin +++)

$\hspace{0.01em}$

{\xssc LABEL: {Remark Infin}}
\label{Remark Infin}

Note that our schemata allow us to generate infintely many new rules, here
is
an example:

Start with A, add $s_{1,1},$ $s_{1,2}$ two sets small in $A \xcv s_{1,1}$
$(A \xcv s_{1,2}$ respectively).
Consider now $A \xcv s_{1,1} \xcv s_{1,2}$ and $s_{2}$ s.t. $s_{2}$ is
small in $A \xcv s_{1,1} \xcv s_{1,2} \xcv s_{2}.$
Continue with $s_{3,1},$ $s_{3,2}$ small in $A \xcv s_{1,1} \xcv s_{1,2}
\xcv s_{2} \xcv s_{3,1}$ etc.

Without additional properties, this system creates a new rule, which is
not equivalent to any usual rules.

$ \xcz $
\\[3ex]
\subsection{
The finite versions
}

\er

\bfa

$\hspace{0.01em}$

(+++ Orig. No.:  Fact I-n +++)

$\hspace{0.01em}$

{\xssc LABEL: {Fact I-n}}
\label{Fact I-n}

(1) $(I_{n})$ $+$ $(eM \xdi )$ $ \xch $ $( \xdm^{+}_{n}),$

(2) $(I_{n})$ $+$ $(eM \xdi )$ $ \xch $ $(CM_{n}),$

(3) $(I_{n})$ $+$ $(eM \xdi )$ $ \xch $ $(OR_{n}).$

\efa

\subparagraph{
Proof
}

$\hspace{0.01em}$

(+++ Orig.:  Proof +++)

(1)

Let $X_{1} \xcc  \Xl  \xcc X_{n},$ so
$X_{n}=X_{1} \xcv (X_{2}-X_{1}) \xcv  \Xl  \xcv (X_{n}-X_{n-1}).$ Let
$X_{i} \xbe \xdf (X_{i+1}),$ so $X_{i+1}-X_{i} \xbe \xdi (X_{i+1}) \xcc
\xdi (X_{n})$
by $(eM \xdi )$ for $1 \xck i \xck n-1,$ so by $(I_{n})$ $X_{1} \xbe
\xdm^{+}(X_{n}).$

(2)

Suppose $ \xba \xcn \xbb_{1}, \Xl, \xba \xcn \xbb_{n-1},$ but $ \xba \xcu
\xbb_{1} \xcu  \Xl  \xcu \xbb_{n-2} \xcn \xCN \xbb_{n-1}.$
Then $M( \xba \xcu \xCN \xbb_{1}), \Xl,M( \xba \xcu \xCN \xbb_{n-1}) \xbe
\xdi (M( \xba )),$ and
$M( \xba \xcu \xbb_{1} \xcu  \Xl  \xcu \xbb_{n-2} \xcu \xbb_{n-1}) \xbe
\xdi (M( \xba \xcu \xbb_{1} \xcu  \Xl  \xcu \xbb_{n-2})) \xcc \xdi (M(
\xba ))$ by $(eM \xdi ).$
But $M( \xba )=M( \xba \xcu \xCN \xbb_{1}) \xcv  \Xl  \xcv M( \xba \xcu
\xCN \xbb_{n-1}) \xcv M( \xba \xcu \xbb_{1} \xcu  \Xl  \xcu \xbb_{n-2}
\xcu \xbb_{n-1})$ is
now the union of $n$ small subsets, $contradiction.$

(3)

Let $ \xba_{1} \xcn \xbb, \Xl, \xba_{n-1} \xcn \xbb,$ so $M( \xba_{i}
\xcu \xCN \xbb ) \xbe \xdi (M( \xba_{i}))$ for $1 \xck i \xck n-1,$ so
$M( \xba_{i} \xcu \xCN \xbb ) \xbe \xdi (M( \xba_{1} \xco  \Xl  \xco
\xba_{n-1}))$ for $1 \xck i \xck n-1$ by $(eM \xdi ),$ so
$M(( \xba_{1} \xco  \Xl  \xco \xba_{n-1}) \xcu \xbb )$ $=$ $M( \xba_{1}
\xco  \Xl  \xco \xba_{n-1})- \xcv \{M( \xba_{i} \xcu \xCN \xbb ):1 \xck i
\xck n-1\}$ $ \xce $
$ \xdi (M( \xba_{1} \xco  \Xl  \xco \xba_{n-1}))$ by $(I_{n}),$ so $
\xba_{1} \xco  \Xl  \xco \xba_{n-1} \xcN \xCN \xbb.$

$ \xcz $
\\[3ex]

In the following example, $(OR_{n}),$ $( \xdm^{+}_{n}),$ $(CM_{n})$ hold,
but $( \xdi_{n})$ fails, so
by Fact \ref{Fact I-n} (page \pageref{Fact I-n})  $( \xdi_{n})$ is strictly
stronger
than $(OR_{n}),$ $( \xdm^{+}_{n}),$ $(CM_{n}).$

\be

$\hspace{0.01em}$

(+++ Orig. No.:  Example Not-I-n +++)

$\hspace{0.01em}$

{\xssc LABEL: {Example Not-I-n}}
\label{Example Not-I-n}

Let $n \xcg 3.$

Consider $X:=\{1, \Xl,n\},$ $ \xdy:= \xdp (X)-\{ \xCQ \},$
$ \xdi (X):=\{ \xCQ \} \xcv \{\{i\}:1 \xck i \xck n\},$ and for all $Y
\xcb X$ $ \xdi (Y):=\{ \xCQ \}.$

$ \xCf (Opt),$ $ \xCf (iM),$ $(eM \xdi ),$ $(eM \xdf )$ (by Fact \ref{Fact
1-element} (page \pageref{Fact 1-element}) ),
$(1*s),$ $(2*s)$ hold, $(I_{n})$ fails, of course.

(1) $(OR_{n})$ holds:

Suppose $ \xba_{1} \xcn \xbb, \Xl, \xba_{n-1} \xcn \xbb,$ $ \xba_{1}
\xco  \Xl  \xco \xba_{n-1} \xcn \xCN \xbb.$

Case 1: $ \xba_{1} \xco  \Xl  \xco \xba_{n-1} \xcl \xCN \xbb,$ then for
all $i$ $ \xba_{i} \xcl \xCN \xbb,$ so for no $i$ $ \xba_{i} \xcn \xbb $
by $(1*s)$ and thus $(AND_{1}),$ $contradiction.$

Case 2: $ \xba_{1} \xco  \Xl  \xco \xba_{n-1} \xcL \xCN \xbb,$ then $M(
\xba_{1} \xco  \Xl  \xco \xba_{n-1})=X,$ and there is exactly
1 $k \xbe X$ s.t. $k \xcm \xbb.$ Fix this $k.$
By prerequisite, $ \xba_{i} \xcn \xbb.$ If $M( \xba_{i})=X,$ $ \xba_{i}
\xcl \xbb $ cannot be, so there must be
exactly 1 $k' $ s.t. $k' \xcm \xCN \xbb,$ but $card(X) \xcg 3,$
$contradiction.$ So $M( \xba_{i}) \xcb X,$ and $ \xba_{i} \xcl \xbb,$ so
$M( \xba_{i})= \xCQ $ or
$M( \xba_{i})=\{k\}$ for all $i,$ so $M( \xba_{1} \xco  \Xl  \xco
\xba_{n-1}) \xEd X,$ $contradiction.$

(2) $( \xdm^{+}_{n})$ holds:

$( \xdm^{+}_{n})$ is a consequence of $( \xdm^{+}_{ \xbo }),$ (3) so it
suffices to show that the latter
holds. Let $X_{1} \xbe \xdf (X_{2}),$ $X_{2} \xbe \xdf (X_{3}).$ Then
$X_{1}=X_{2}$ or $X_{2}=X_{3},$ so the result is
trivial.

(3) $(CM_{n})$ holds:

Suppose $ \xba \xcn \xbb_{1}, \Xl, \xba \xcn \xbb_{n-1},$ $ \xba \xcu
\xbb_{1} \xcu  \Xl  \xcu \xbb_{n-2} \xcn \xCN \xbb_{n-1}.$

Case 1: For all $i,$ $1 \xck i \xck n-2,$ $ \xba \xcl \xbb_{i},$ then $M(
\xba \xcu \xbb_{1} \xcu  \Xl  \xcu \xbb_{n-2})=M( \xba ),$
so $ \xba \xcn \xbb_{n-1}$ and $ \xba \xcn \xCN \xbb_{n-1},$
$contradiction.$

Case 2: There is $i,$ $1 \xck i \xck n-2,$ $ \xba \xcL \xbb_{i},$ then $M(
\xba )=X,$
$M( \xba \xcu \xbb_{1} \xcu  \Xl  \xcu \xbb_{n-2}) \xcb M( \xba ),$ so $
\xba \xcu \xbb_{1} \xcu  \Xl  \xcu \xbb_{n-2} \xcl \xCN \xbb_{n-1}.$
$Card(M( \xba \xcu \xbb_{1} \xcu  \Xl  \xcu \xbb_{n-2})) \xcg n-(n-2)=2,$
so $card(M( \xCN \xbb_{n-1})) \xcg 2,$ so
$ \xba \xcN \xbb_{n-1},$ $contradiction.$

$ \xcz $
\\[3ex]
\subsection{
The $\xbo$ version
}

\ee

\bfa

$\hspace{0.01em}$

(+++ Orig. No.:  Fact CM-Omega +++)

$\hspace{0.01em}$

{\xssc LABEL: {Fact CM-Omega}}
\label{Fact CM-Omega}

$(CM_{ \xbo })$ $ \xcj $ $( \xdm^{+}_{ \xbo })$ (4)

\efa

\subparagraph{
Proof
}

$\hspace{0.01em}$

(+++ Orig.:  Proof +++)

`` $ \xch $ ''

Suppose all sets are definable.

Let $A,B \xbe \xdi (X),$
$X=M( \xba ),$ $A=M( \xba \xcu \xCN \xbb ),$ $B=M( \xba \xcu \xCN \xbb '
),$ so $ \xba \xcn \xbb,$ $ \xba \xcn \xbb ',$ so by $(CM_{ \xbo })$
$ \xba \xcu \xbb ' \xcn \xbb,$ so $A-B=M( \xba \xcu \xbb ' \xcu \xCN \xbb
) \xbe \xdi (M( \xba \xcu \xbb ' ))= \xdi (X-$B).

`` $ \xci $ ''

Let $ \xba \xcn \xbb,$ $ \xba \xcn \xbb ',$ so $M( \xba \xcu \xCN \xbb )
\xbe \xdi (M( \xba )),$ $M( \xba \xcu \xCN \xbb ' ) \xbe \xdi (M( \xba
)),$ so
by prerequisite $M( \xba \xcu \xCN \xbb ' )-M( \xba \xcu \xCN \xbb )=M(
\xba \xcu \xbb \xcu \xCN \xbb ' )$ $ \xbe $
$ \xdi (M( \xba )-M( \xba \xcu \xCN \xbb ))= \xdi (M( \xba \xcu \xbb )),$
so $ \xba \xcu \xbb \xcn \xbb '.$

$ \xcz $
\\[3ex]

\bfa

$\hspace{0.01em}$

(+++ Orig. No.:  Fact I-Omega +++)

$\hspace{0.01em}$

{\xssc LABEL: {Fact I-Omega}}
\label{Fact I-Omega}

(1) $(I_{ \xbo })$ $+$ $(eM \xdi )$ $ \xch $ $(OR_{ \xbo }),$

(2) $(I_{ \xbo })$ $+$ $(eM \xdi )$ $ \xch $ $( \xdm^{+}_{ \xbo })$ (1),

(3) $(I_{ \xbo })$ $+$ $(eM \xdf )$ $ \xch $ $( \xdm^{+}_{ \xbo })$ (2),

(4) $(I_{ \xbo })$ $+$ $(eM \xdi )$ $ \xch $ $( \xdm^{+}_{ \xbo })$ (3),

(5) $(I_{ \xbo })$ $+$ $(eM \xdf )$ $ \xch $ $( \xdm^{+}_{ \xbo })$ (4)
(and thus, by
Fact \ref{Fact CM-Omega} (page \pageref{Fact CM-Omega}), $(CM_{ \xbo })).$

\efa

\subparagraph{
Proof
}

$\hspace{0.01em}$

(+++ Orig.:  Proof +++)

(1)

Let $ \xba \xcn \xbb,$ $ \xba ' \xcn \xbb $ $ \xch $ $M( \xba \xcu \xCN
\xbb ) \xbe \xdi (M( \xba )),$ $M( \xba ' \xcu \xCN \xbb ) \xbe \xdi (M(
\xba ' )),$
so by $(eM \xdi )$ $M( \xba \xcu \xCN \xbb ) \xbe \xdi (M( \xba \xco \xba
' )),$ $M( \xba ' \xcu \xCN \xbb ) \xbe \xdi (M( \xba \xco \xba ' )),$
so $M(( \xba \xco \xba ' ) \xcu \xCN \xbb ) \xbe \xdi (M( \xba \xco \xba '
))$ by $(I_{ \xbo }),$ so $ \xba \xco \xba ' \xcn \xbb.$

(2)

Let $A \xcc X \xcc Y,$
$A \xbe \xdi (Y),$ $X-A \xbe \xdi (X) \xcc_{(eM \xdi )} \xdi (Y)$ $ \xch $
$X=(X-A) \xcv A \xbe \xdi (Y)$ by $(I_{ \xbo }).$

(3)

Let $A \xcc X \xcc Y,$
let $A \xbe \xdi (Y),$ $Y-X \xbe \xdi (Y)$ $ \xch $ $A \xcv (Y-X) \xbe
\xdi (Y)$ by $(I_{ \xbo })$ $ \xch $
$X-A=Y-(A \xcv (Y-X)) \xbe \xdf (Y)$ $ \xch $ $X-A \xbe \xdf (X)$ by $(eM
\xdf ).$

(4)

Let $A \xcc X \xcc Y,$ $A \xbe \xdf (X),$ $X \xbe \xdf (Y),$ so
$Y-X \xbe \xdi (Y),$ $X-A \xbe \xdi (X) \xcc_{(eM \xdi )} \xdi (Y)$ $ \xch
$ $Y-A=(Y-X) \xcv (X-A) \xbe \xdi (Y)$ by $( \xdi_{ \xbo })$ $ \xch $
$A \xbe \xdf (Y).$

(5)

Let $A,B \xcc X,$
$A,B \xbe \xdi (X)$ $ \xch_{(I_{ \xbo })}$ $A \xcv B \xbe \xdi (X)$ $ \xch
$ $X-(A \xcv B) \xbe \xdf (X),$ but $X-(A \xcv B) \xcc X-$B,
so $X-(A \xcv B) \xbe \xdf (X-$B) by $(eM \xdf ),$ so $A-B=(X-B)-(X-(A
\xcv B)) \xbe \xdi (X-$B).

$ \xcz $
\\[3ex]

We give three examples of independence of the various versions of
$( \xdm^{+}_{ \xbo }).$

\be

$\hspace{0.01em}$

(+++ Orig. No.:  Example Versions-M-Omega +++)

$\hspace{0.01em}$

{\xssc LABEL: {Example Versions-M-Omega}}
\label{Example Versions-M-Omega}

All numbers refer to the versions of $( \xdm^{+}_{ \xbo }).$

For easier reading, we re-write for $A \xcc X \xcc Y$

$( \xdm^{+}_{ \xbo })(1):$ $A \xbe \xdf (X),$ $A \xbe \xdi (Y)$ $ \xch $
$X \xbe \xdi (Y),$

$( \xdm^{+}_{ \xbo })(2):$ $X \xbe \xdf (Y),$ $A \xbe \xdi (Y)$ $ \xch $
$A \xbe \xdi (X).$

We give three examples. Investigating all possibilities exhaustively
seems quite tedious, and might best be done with the help of a computer.
Fact \ref{Fact 1-element} (page \pageref{Fact 1-element})  will be used
repeatedly.

 \xEI

 \xDH

(1), (2), (4) fail, (3) holds:

Let $Y:=\{a,b,c\},$ $ \xdy:= \xdp (Y)-\{ \xCQ \},$ $ \xdf
(Y):=\{\{a,c\},$ $\{b,c\},$ $Y\}$

Let $X:=\{a,b\},$ $ \xdf (X):=\{\{a\},$ $X\},$ $A:=\{a\},$ and $ \xdf
(Z):=\{Z\}$ for all $Z \xEd X,Y.$

$ \xCf (Opt),$ $ \xCf (iM),$ $(eM \xdi ),$ $(eM \xdf )$ hold, $(I_{ \xbo
})$ fails, of course.

(1) fails: $A \xbe \xdf (X),$ $A \xbe \xdi (Y),$ $X \xce \xdi (Y).$

(2) fails: $\{a,c\} \xbe \xdf (Y),$ $\{a\} \xbe \xdi (Y),$ but $\{a\} \xce
\xdi (\{a,c\}).$

(3) holds: If $X_{1} \xbe \xdf (X_{2}),$ $X_{2} \xbe \xdf (X_{3}),$ then
$X_{1}=X_{2}$ or $X_{2}=X_{3},$ so (3) holds
trivially (note that $X \xce \xdf (Y)).$

(4) fails: $\{a\},\{b\} \xbe \xdi (Y),$ $\{a\} \xce \xdi (Y-\{b\})= \xdi
(\{a,c\})=\{ \xCQ \}.$

 \xDH

(2), (3), (4) fail, (1) holds:

Let $Y:=\{a,b,c\},$ $ \xdy:= \xdp (Y)-\{ \xCQ \},$ $ \xdf
(Y):=\{\{a,b\},$ $\{a,c\},$ $Y\}$

Let $X:=\{a,b\},$ $ \xdf (X):=\{\{a\},$ $X\},$ and $ \xdf (Z):=\{Z\}$ for
all $Z \xEd X,Y.$

$ \xCf (Opt),$ $ \xCf (iM),$ $(eM \xdi ),$ $(eM \xdf )$ hold, $(I_{ \xbo
})$ fails, of course.

(1) holds:

Let $X_{1} \xbe \xdf (X_{2}),$ $X_{1} \xbe \xdi (X_{3}),$ we have to show
$X_{2} \xbe \xdi (X_{3}).$
If $X_{1}=X_{2},$ then this is trivial. Consider $X_{1} \xbe \xdf
(X_{2}).$
If $X_{1} \xEd X_{2},$ then $X_{1}$ has to be $\{a\}$ or
$\{a,b\}$ or $\{a,c\}.$ But none of these are in $ \xdi (X_{3})$ for any
$X_{3},$ so the
implication is trivially true.

(2) fails: $\{a,c\} \xbe \xdf (Y),$ $\{c\} \xbe \xdi (Y),$ $\{c\} \xce
\xdi (\{a,c\}).$

(3) fails: $\{a\} \xbe \xdf (X),$ $X \xbe \xdf (Y),$ $\{a\} \xce \xdf
(Y).$

(4) fails: $\{b\},\{c\} \xbe \xdi (Y),$ $\{c\} \xce \xdi (Y-\{b\})= \xdi
(\{a,c\})=\{ \xCQ \}.$

 \xDH

(1), (2), (4) hold, (3) fails:

Let $Y:=\{a,b,c\},$ $ \xdy:= \xdp (Y)-\{ \xCQ \},$ $ \xdf
(Y):=\{\{a,b\},$ $\{a,c\},$ $Y\}$

Let $ \xdf (\{a,b\}):=\{\{a\},\{a,b\}\},$ $ \xdf
(\{a,c\}):=\{\{a\},\{a,c\}\},$
and $ \xdf (Z):=\{Z\}$ for all other $Z.$

$ \xCf (Opt),$ $ \xCf (iM),$ $(eM \xdi ),$ $(eM \xdf )$ hold, $(I_{ \xbo
})$ fails, of course.

(1) holds:

Let $X_{1} \xbe \xdf (X_{2}),$ $X_{1} \xbe \xdi (X_{3}),$ we have to show
$X_{2} \xbe \xdi (X_{3}).$ Consider $X_{1} \xbe \xdi (X_{3}).$
If $X_{1}=X_{2},$ this is trivial. If $ \xCQ \xEd X_{1} \xbe \xdi
(X_{3}),$ then $X_{1}=\{b\}$ or $X_{1}=\{c\},$ but
then by $X_{1} \xbe \xdf (X_{2})$ $X_{2}$ has to be $\{b\},$ or $\{c\},$
so $X_{1}=X_{2}.$

(2) holds:
Let $X_{1} \xcc X_{2} \xcc X_{3},$ let $X_{2} \xbe \xdf (X_{3}),$ $X_{1}
\xbe \xdi (X_{3}),$ we have to show $X_{1} \xbe \xdi (X_{2}).$
If $X_{1}= \xCQ,$ this is trivial, likewise if $X_{2}=X_{3}.$ Otherwise
$X_{1}=\{b\}$ or $X_{1}=\{c\},$ and $X_{3}=Y.$ If $X_{1}=\{b\},$ then
$X_{2}=\{a,b\},$ and the condition
holds, likewise if $X_{1}=\{c\},$ then $X_{2}=\{a,c\},$ and it holds
again.

(3) fails: $\{a\} \xbe \xdf (\{a,c\}),$ $\{a,c\} \xbe \xdf (Y),$ $\{a\}
\xce \xdf (Y).$

(4) holds:

If $A,B \xbe \xdi (X),$ and $A \xEd B,$ $A,B \xEd \xCQ,$ then $X=Y$ and
e.g. $A=\{c\},$ $B=\{b\},$ and
$\{c\} \xbe \xdi (Y-\{b\})= \xdi (\{a,c\}).$

 \xEJ
$ \xcz $
\\[3ex]
\subsection{
Rational Monotony
}

\ee

\bfa

$\hspace{0.01em}$

(+++ Orig. No.:  Fact M-plus-plus +++)

$\hspace{0.01em}$

{\xssc LABEL: {Fact M-plus-plus}}
\label{Fact M-plus-plus}

The three versions of $( \xdm^{++})$ are equivalent.

(We assume closure of the domain under set difference.
For the third version of $( \xdm^{++}),$ we use $ \xCf (iM).)$

\efa

\subparagraph{
Proof
}

$\hspace{0.01em}$

(+++ Orig.:  Proof +++)

For (1) and (2), we have $A,B \xcc X,$ for (3) we have $A \xcc X \xcc Y.$
For $A,B \xcc X,$ $(X-B)-((X-A)-B)=A-B$ holds.

$(1) \xch (2):$ Let $A \xbe \xdf (X),$ $B \xce \xdf (X),$ so $X-A \xbe
\xdi (X),$ so by prerequisite
$(X-A)-B \xbe \xdi (X-$B), so $A-B=(X-B)-((X-A)-B) \xbe \xdf (X-$B).

$(2) \xch (1):$ Let $A \xbe \xdi (X),$ $B \xce \xdf (X),$ so $X-A \xbe
\xdf (X),$ so by prerequisite
$(X-A)-B \xbe \xdf (X-$B), so $A-B=(X-B)-((X-A)-B) \xbe \xdi (X-$B).

$(1) \xch (3):$

Suppose $A \xce \xdm^{+}(Y),$ but $X \xbe \xdm^{+}(Y),$ we show $A \xce
\xdm^{+}(X).$ So $A \xbe \xdi (Y),$ $Y-X \xce \xdf (Y),$
so by (1) $A=A-(Y-X) \xbe \xdi (Y-(Y-X))= \xdi (X).$

$(3) \xch (1):$

Suppose $A-B \xce \xdi (X-$B), $B \xce \xdf (X),$ we show $A \xce \xdi
(X).$ By prerequisite $A-B \xbe \xdm^{+}(X-$B),
$X-B \xbe \xdm^{+}(X),$ so by (3) $A-B \xbe \xdm^{+}(X),$ so by $ \xCf
(iM)$ $A \xbe \xdm^{+}(X),$ so $A \xce \xdi (X).$

$ \xcz $
\\[3ex]

\bfa

$\hspace{0.01em}$

(+++ Orig. No.:  Fact M-RatM +++)

$\hspace{0.01em}$

{\xssc LABEL: {Fact M-RatM}}
\label{Fact M-RatM}

We assume that all sets are definable by a formula.

$ \xCf (RatM)$ $ \xcj $ $( \xdm^{++})$

\efa

\subparagraph{
Proof
}

$\hspace{0.01em}$

(+++ Orig.:  Proof +++)

We show equivalence of $ \xCf (RatM)$ with version (1) of $( \xdm^{++}).$

`` $ \xch $ ''

We have $A,B \xcc X,$ so we can write
$X=M( \xbf ),$ $A=M( \xbf \xcu \xCN \xbq ),$ $B=M( \xbf \xcu \xCN \xbq '
).$ $A \xbe \xdi (X),$ $B \xce \xdf (X),$ so
$ \xbf \xcn \xbq,$ $ \xbf \xcN \xCN \xbq ',$ so by $ \xCf (RatM)$ $ \xbf
\xcu \xbq ' \xcn \xbq,$ so
$A-B=M( \xbf \xcu \xCN \xbq )-M( \xbf \xcu \xCN \xbq ' )=M( \xbf \xcu \xbq
' \xcu \xCN \xbq ) \xbe \xdi (M( \xbf \xcu \xbq ' ))= \xdi (X-$B).

`` $ \xci $ ''

Let $ \xbf \xcn \xbq,$ $ \xbf \xcN \xCN \xbq ',$ so $M( \xbf \xcu \xCN
\xbq ) \xbe \xdi (M( \xbf )),$ $M( \xbf \xcu \xCN \xbq ' ) \xce \xdf (M(
\xbf )),$ so
by $( \xdm^{++})$ (1) $M( \xbf \xcu \xbq ' \xcu \xCN \xbq )=M( \xbf \xcu
\xCN \xbq )-M( \xbf \xcu \xCN \xbq ' ) \xbe \xdi (M( \xbf \xcu \xbq ' )),$
so
$ \xbf \xcu \xbq ' \xcn \xbq.$

$ \xcz $
\\[3ex]
\section{
Size and principal filter logic
}

$\hspace{0.01em}$

{\xssc LABEL: {Section Principal}}
\label{Section Principal}

The connection with logical rules is shown in the following table
Definition \ref{Definition Log-Cond-Ref-Size} (page \pageref{Definition
Log-Cond-Ref-Size}).
Most of the table was already published in  \cite{GS08c},
it is repeated here for the reader's convenience.
\index{Definition Log-Cond-Ref-Size}

$\hspace{0.01em}$

{\xssc LABEL: {Definition Log-Cond-Ref-Size}}
\label{Definition Log-Cond-Ref-Size}

The numbers in the first column ``Correspondence''
refer to
Proposition 21 in  \cite{GS08c},
those in the second column ``Correspondence''
to Proposition \ref{Proposition Ref-Class-Mu-neu} (page \pageref{Proposition
Ref-Class-Mu-neu}).

\begin{turn}{90}


{\xssc

\begin{tabular}{|c|c|c|c|c|c|}

\hline

\multicolumn{2}{|c|}{Logical rule}
\xEH
Correspondence
\xEH
Model set
\xEH
Correspondence
\xEH
Size Rules
\xEP

\hline

\multicolumn{6}{|c|}{Basics}
\xEP

\hline

$(SC)$ Supraclassicality
\xEH
$(SC)$
\xEH
$\xch$ (4.1)
\xEH
$( \xbm \xcc )$
\xEH
trivial
\xEH
$(Opt)$
\xEP

\cline{3-3}

$ \xbf \xcl \xbq $ $ \xch $ $ \xbf \xcn \xbq $
\xEH
$ \ol{T} \xcc \ol{ \ol{T} }$
\xEH
$\xci$ (4.2)
\xEH
$f(X) \xcc X$
\xEH
\xEH
\xEP

\cline{1-1}

$(REF)$ Reflexivity
\xEH
\xEH
\xEH
\xEH
\xEH
\xEP

$ T \xcv \{\xba\} \xcn \xba $
\xEH
\xEH
\xEH
\xEH
\xEH
\xEP

\hline

$(LLE)$
\xEH
$(LLE)$
\xEH
\xEH
\xEH
\xEH
\xEP

Left Logical Equivalence
\xEH
\xEH
\xEH
\xEH
\xEH
\xEP

$ \xcl \xbf \xcr \xbf ',  \xbf \xcn \xbq   \xch $
\xEH
$ \ol{T}= \ol{T' }  \xch   \ol{\ol{T}} = \ol{\ol{T'}}$
\xEH
\xEH
\xEH
\xEH
\xEP

$ \xbf ' \xcn \xbq $
\xEH
\xEH
\xEH
\xEH
\xEH
\xEP

\hline

$(RW)$ Right Weakening
\xEH
$(RW)$
\xEH
\xEH
\xEH
trivial
\xEH
$(iM)$
\xEP

$ \xbf \xcn \xbq,  \xcl \xbq \xcp \xbq '   \xch $
\xEH
$ T \xcn \xbq,  \xcl \xbq \xcp \xbq '   \xch $
\xEH
\xEH
\xEH
\xEH
\xEP

$ \xbf \xcn \xbq ' $
\xEH
$T \xcn \xbq ' $
\xEH
\xEH
\xEH
\xEH
\xEP

\hline

$(wOR)$
\xEH
$(wOR)$
\xEH
$\xch$ (3.1)
\xEH
$( \xbm wOR)$
\xEH
$\xci$ (1.1)
\xEH
$(eM\xdi)$
\xEP

\cline{3-3}
\cline{5-5}

$ \xbf \xcn \xbq,$ $ \xbf ' \xcl \xbq $ $ \xch $
\xEH
$ \ol{ \ol{T} } \xcs \ol{T' }$ $ \xcc $ $ \ol{ \ol{T \xco T' } }$
\xEH
$\xci$ (3.2)
\xEH
$f(X \xcv Y) \xcc f(X) \xcv Y$
\xEH
$\xch$ (1.2)
\xEH
\xEP

$ \xbf \xco \xbf ' \xcn \xbq $
\xEH
\xEH
\xEH
\xEH
\xEH
\xEP

\hline

$(disjOR)$
\xEH
$(disjOR)$
\xEH
$\xch$ (2.1)
\xEH
$( \xbm disjOR)$
\xEH
$\xci$ (4.1)
\xEH
$(I\xcv disj)$
\xEP

\cline{3-3}
\cline{5-5}

$ \xbf \xcl \xCN \xbf ',$ $ \xbf \xcn \xbq,$
\xEH
$\xCN Con(T \xcv T') \xch$
\xEH
$\xci$ (2.2)
\xEH
$X \xcs Y= \xCQ $ $ \xch $
\xEH
$\xch$ (4.2)
\xEH
\xEP

$ \xbf ' \xcn \xbq $ $ \xch $ $ \xbf \xco \xbf ' \xcn \xbq $
\xEH
$ \ol{ \ol{T} } \xcs \ol{ \ol{T' } } \xcc \ol{ \ol{T \xco T' } }$
\xEH
\xEH
$f(X \xcv Y) \xcc f(X) \xcv f(Y)$
\xEH
\xEH
\xEP

\hline

$(CP)$
\xEH
$(CP)$
\xEH
$\xch$ (5.1)
\xEH
$( \xbm \xCQ )$
\xEH
trivial
\xEH
$(I_1)$
\xEP

\cline{3-3}

Consistency Preservation
\xEH
\xEH
$\xci$ (5.2)
\xEH
\xEH
\xEH
\xEP

$ \xbf \xcn \xcT $ $ \xch $ $ \xbf \xcl \xcT $
\xEH
$T \xcn \xcT $ $ \xch $ $T \xcl \xcT $
\xEH
\xEH
$f(X)= \xCQ $ $ \xch $ $X= \xCQ $
\xEH
\xEH
\xEP

\hline

\xEH
\xEH
\xEH
$( \xbm \xCQ fin)$
\xEH
\xEH
$(I_1)$
\xEP

\xEH
\xEH
\xEH
$X \xEd \xCQ $ $ \xch $ $f(X) \xEd \xCQ $
\xEH
\xEH
\xEP

\xEH
\xEH
\xEH
for finite $X$
\xEH
\xEH
\xEP

\hline

\xEH
$(AND_1)$
\xEH
\xEH
\xEH
\xEH
$(I_2)$
\xEP

\xEH
$\xba\xcn\xbb \xch \xba\xcN\xCN\xbb$
\xEH
\xEH
\xEH
\xEH
\xEP

\hline

\xEH
$(AND_n)$
\xEH
\xEH
\xEH
\xEH
$(I_n)$
\xEP

\xEH
$\xba\xcn\xbb_1, \ldots, \xba\xcn\xbb_{n-1} \xch $
\xEH
\xEH
\xEH
\xEH
\xEP

\xEH
$\xba\xcN(\xCN\xbb_1 \xco \ldots \xco \xCN\xbb_{n-1})$
\xEH
\xEH
\xEH
\xEH
\xEP

\hline

$(AND)$
\xEH
$(AND)$
\xEH
\xEH
\xEH
trivial
\xEH
$(I_\xbo)$
\xEP

$ \xbf \xcn \xbq,  \xbf \xcn \xbq '   \xch $
\xEH
$ T \xcn \xbq, T \xcn \xbq '   \xch $
\xEH
\xEH
\xEH
\xEH
\xEP

$ \xbf \xcn \xbq \xcu \xbq ' $
\xEH
$ T \xcn \xbq \xcu \xbq ' $
\xEH
\xEH
\xEH
\xEH
\xEP

\hline

$(CCL)$ Classical Closure
\xEH
$(CCL)$
\xEH
\xEH
\xEH
trivial
\xEH
$(iM)+(I_\xbo)$
\xEP

\xEH
$ \ol{ \ol{T} }$ classically closed
\xEH
\xEH
\xEH
\xEH
\xEP

\hline

$(OR)$
\xEH
$(OR)$
\xEH
$\xch$ (1.1)
\xEH
$( \xbm OR)$
\xEH
$\xci$ (2.1)
\xEH
$(eM\xdi)+(I_\xbo)$
\xEP

\cline{3-3}
\cline{5-5}

$ \xbf \xcn \xbq,  \xbf ' \xcn \xbq   \xch $
\xEH
$ \ol{\ol{T}} \xcs \ol{\ol{T'}} \xcc \ol{\ol{T \xco T'}} $
\xEH
$\xci$ (1.2)
\xEH
$f(X \xcv Y) \xcc f(X) \xcv f(Y)$
\xEH
$\xch$ (2.2)
\xEH
\xEP

$ \xbf \xco \xbf ' \xcn \xbq $
\xEH
\xEH
\xEH
\xEH
\xEH
\xEP

\hline

\xEH
$(PR)$
\xEH
$\xch$ (6.1)
\xEH
$( \xbm PR)$
\xEH
$\xci$ (3.1)
\xEH
$(eM\xdi)+(I_\xbo)$
\xEP

\cline{3-3}
\cline{5-5}

$ \ol{ \ol{ \xbf \xcu \xbf ' } }$ $ \xcc $ $ \ol{ \ol{ \ol{ \xbf } } \xcv
\{ \xbf ' \}}$
\xEH
$ \ol{ \ol{T \xcv T' } }$ $ \xcc $ $ \ol{ \ol{ \ol{T} } \xcv T' }$
\xEH
$\xci (\xbm dp)+(\xbm\xcc)$ (6.2)
\xEH
$X \xcc Y$ $ \xch $
\xEH
$\xch$ (3.2)
\xEH
\xEP

\cline{3-3}

\xEH
\xEH
$\xcI$ without $(\xbm dp)$ (6.3)
\xEH
$f(Y) \xcs X \xcc f(X)$
\xEH
\xEH
\xEP

\cline{3-3}

\xEH
\xEH
$\xci (\xbm\xcc)$ (6.4)
\xEH
\xEH
\xEH
\xEP

\xEH
\xEH
$T'$ a formula
\xEH
\xEH
\xEH
\xEP

\cline{3-4}

\xEH
\xEH
$\xci$ (6.5)
\xEH
$(\xbm PR ')$
\xEH
\xEH
\xEP

\xEH
\xEH
$T'$ a formula
\xEH
$f(X) \xcs Y \xcc f(X \xcs Y)$
\xEH
\xEH
\xEP

\hline

$(CUT)$
\xEH
$(CUT)$
\xEH
$\xch$ (7.1)
\xEH
$ (\xbm CUT) $
\xEH
$\xci$ (8.1)
\xEH
$(eM\xdi)+(I_\xbo)$
\xEP

\cline{3-3}
\cline{5-5}

$ T  \xcn \xba; T \xcv \{ \xba\} \xcn \xbb \xch $
\xEH
$T \xcc \ol{T' } \xcc \ol{ \ol{T} }  \xch $
\xEH
$\xci$ (7.2)
\xEH
$f(X) \xcc Y \xcc X  \xch $
\xEH
$\xcH$ (8.2)
\xEH
\xEP

$ T  \xcn \xbb $
\xEH
$ \ol{ \ol{T'} } \xcc \ol{ \ol{T} }$
\xEH
\xEH
$f(X) \xcc f(Y)$
\xEH
\xEH
\xEP

\hline

\end{tabular}

}

\end{turn}

\begin{turn}{90}

{\xssc

\begin{tabular}{|c|c|c|c|c|c|}

\hline

\multicolumn{2}{|c|}{Logical rule}
\xEH
Correspondence
\xEH
Model set
\xEH
Correspondence
\xEH
Size-Rule
\xEP

\hline

\multicolumn{6}{|c|}{Cumulativity}
\xEP

\hline

$(wCM)$
\xEH
\xEH
\xEH
\xEH
trivial
\xEH
$(eM\xdf)$
\xEP

$\xba\xcn\xbb, \xba'\xcl\xba, \xba\xcu\xbb\xcl\xba' \xch \xba'\xcn\xbb$
\xEH
\xEH
\xEH
\xEH
\xEH
\xEP

\hline

$(CM_2)$
\xEH
\xEH
\xEH
\xEH
\xEH
$(I_2)$
\xEP

$\xba\xcn\xbb, \xba\xcn\xbb' \xch \xba\xcu\xbb\xcL\xCN\xbb'$
\xEH
\xEH
\xEH
\xEH
\xEH
\xEP

\hline

$(CM_n)$
\xEH
\xEH
\xEH
\xEH
\xEH
$(I_n)$
\xEP

$\xba\xcn\xbb_1, \ldots, \xba\xcn\xbb_n \xch $
\xEH
\xEH
\xEH
\xEH
\xEH
\xEP

$\xba \xcu \xbb_1 \xcu \ldots \xcu \xbb_{n-1} \xcL\xCN\xbb_n$
\xEH
\xEH
\xEH
\xEH
\xEH
\xEP

\hline

$(CM)$ Cautious Monotony
\xEH
$(CM)$
\xEH
$\xch$ (8.1)
\xEH
$ (\xbm CM) $
\xEH
$\xci$ (5.1)
\xEH
$(\xdm^+_\xbo)(4)$
\xEP

\cline{3-3}
\cline{5-5}

$ \xbf \xcn \xbq,  \xbf \xcn \xbq '   \xch $
\xEH
$T \xcc \ol{T' } \xcc \ol{ \ol{T} }  \xch $
\xEH
$\xci$ (8.2)
\xEH
$f(X) \xcc Y \xcc X  \xch $
\xEH
$\xch$ (5.2)
\xEH
\xEP

$ \xbf \xcu \xbq \xcn \xbq ' $
\xEH
$ \ol{ \ol{T} } \xcc \ol{ \ol{T' } }$
\xEH
\xEH
$f(Y) \xcc f(X)$
\xEH
\xEH
\xEP

\cline{1-1}

\cline{3-4}

or $(ResM)$ Restricted Monotony
\xEH
\xEH
$\xch$ (9.1)
\xEH
$(\xbm ResM)$
\xEH
\xEH
\xEP

\cline{3-3}

$ T  \xcn \xba, \xbb \xch T \xcv \{\xba\} \xcn \xbb $
\xEH
\xEH
$\xci$ (9.2)
\xEH
$ f(X) \xcc A \xcs B \xch f(X \xcs A) \xcc B $
\xEH
\xEH
\xEP

\hline

$(CUM)$ Cumulativity
\xEH
$(CUM)$
\xEH
$\xch$ (11.1)
\xEH
$( \xbm CUM)$
\xEH
$\xci$ (9.1)
\xEH
$(eM\xdi)+(I_\xbo)+(\xdm^{+}_{\xbo})(4)$
\xEP

\cline{3-3}
\cline{5-5}

$ \xbf \xcn \xbq   \xch $
\xEH
$T \xcc \ol{T' } \xcc \ol{ \ol{T} }  \xch $
\xEH
$\xci$ (11.2)
\xEH
$f(X) \xcc Y \xcc X  \xch $
\xEH
$\xcH$ (9.2)
\xEH
\xEP

$( \xbf \xcn \xbq '   \xcj   \xbf \xcu \xbq \xcn \xbq ' )$
\xEH
$ \ol{ \ol{T} }= \ol{ \ol{T' } }$
\xEH
\xEH
$f(Y)=f(X)$
\xEH
\xEH
\xEP

\hline

\xEH
$ (\xcc \xcd) $
\xEH
$\xch$ (10.1)
\xEH
$ (\xbm \xcc \xcd) $
\xEH
$\xci$ (10.1)
\xEH
$(eM\xdi)+(I_\xbo)+(eM\xdf)$
\xEP

\cline{3-3}
\cline{5-5}

\xEH
$T \xcc \ol{\ol{T'}}, T' \xcc \ol{\ol{T}} \xch $
\xEH
$\xci$ (10.2)
\xEH
$ f(X) \xcc Y, f(Y) \xcc X \xch $
\xEH
$\xcH$ (10.2)
\xEH
\xEP

\xEH
$ \ol{\ol{T'}} = \ol{\ol{T}}$
\xEH
\xEH
$ f(X)=f(Y) $
\xEH
\xEH
\xEP

\hline

\multicolumn{6}{|c|}{Rationality}
\xEP

\hline

$(RatM)$ Rational Monotony
\xEH
$(RatM)$
\xEH
$\xch$ (12.1)
\xEH
$( \xbm RatM)$
\xEH
$\xci$ (6.1)
\xEH
$(\xdm^{++})$
\xEP

\cline{3-3}
\cline{5-5}

$ \xbf \xcn \xbq,  \xbf \xcN \xCN \xbq '   \xch $
\xEH
$Con(T \xcv \ol{\ol{T'}})$, $T \xcl T'$ $ \xch $
\xEH
$\xci$ $(\xbm dp)$ (12.2)
\xEH
$X \xcc Y, X \xcs f(Y) \xEd \xCQ   \xch $
\xEH
$\xch$ (6.2)
\xEH
\xEP

\cline{3-3}

$ \xbf \xcu \xbq ' \xcn \xbq $
\xEH
$ \ol{\ol{T}} \xcd \ol{\ol{\ol{T'}} \xcv T} $
\xEH
$\xcI$ without $(\xbm dp)$ (12.3)
\xEH
$f(X) \xcc f(Y) \xcs X$
\xEH
\xEH
\xEP

\cline{3-3}

\xEH
\xEH
$\xci$ $T$ a formula (12.4)
\xEH
\xEH
\xEH
\xEP

\hline

\xEH
$(RatM=)$
\xEH
$\xch$ (13.1)
\xEH
$( \xbm =)$
\xEH
\xEH
\xEP

\cline{3-3}

\xEH
$Con(T \xcv \ol{\ol{T'}})$, $T \xcl T'$ $ \xch $
\xEH
$\xci$ $(\xbm dp)$ (13.2)
\xEH
$X \xcc Y, X \xcs f(Y) \xEd \xCQ   \xch $
\xEH
\xEH
\xEP

\cline{3-3}

\xEH
$ \ol{\ol{T}} = \ol{\ol{\ol{T'}} \xcv T} $
\xEH
$\xcI$ without $(\xbm dp)$ (13.3)
\xEH
$f(X) = f(Y) \xcs X$
\xEH
\xEH
\xEP

\cline{3-3}

\xEH
\xEH
$\xci$ $T$ a formula (13.4)
\xEH
\xEH
\xEH
\xEP

\hline

\xEH
$(Log=' )$
\xEH
$\xch$ (14.1)
\xEH
$( \xbm =' )$
\xEH
\xEH
\xEP

\cline{3-3}

\xEH
$Con( \ol{ \ol{T' } } \xcv T)$ $ \xch $
\xEH
$\xci$ $(\xbm dp)$ (14.2)
\xEH
$f(Y) \xcs X \xEd \xCQ $ $ \xch $
\xEH
\xEH
\xEP

\cline{3-3}

\xEH
$ \ol{ \ol{T \xcv T' } }= \ol{ \ol{ \ol{T' } } \xcv T}$
\xEH
$\xcI$ without $(\xbm dp)$ (14.3)
\xEH
$f(Y \xcs X)=f(Y) \xcs X$
\xEH
\xEH
\xEP

\cline{3-3}

\xEH
\xEH
$\xci$ $T$ a formula (14.4)
\xEH
\xEH
\xEH
\xEP

\hline

\xEH
$(Log \xFO )$
\xEH
$\xch$ (15.1)
\xEH
$( \xbm \xFO )$
\xEH
\xEH
\xEP

\cline{3-3}

\xEH
$ \ol{ \ol{T \xco T' } }$ is one of
\xEH
$\xci$ (15.2)
\xEH
$f(X \xcv Y)$ is one of
\xEH
\xEH
\xEP

\xEH
$\ol{\ol{T}},$ or $\ol{\ol{T'}},$ or $\ol{\ol{T}} \xcs \ol{\ol{T'}}$ (by (CCL))
\xEH
\xEH
$f(X),$ $f(Y)$ or $f(X) \xcv f(Y)$
\xEH
\xEH
\xEP

\hline

\xEH
$(Log \xcv )$
\xEH
$\xch$ $(\xbm\xcc)+(\xbm=)$ (16.1)
\xEH
$( \xbm \xcv )$
\xEH
\xEH
\xEP

\cline{3-3}

\xEH
$Con( \ol{ \ol{T' } } \xcv T),$ $ \xCN Con( \ol{ \ol{T' } }
\xcv \ol{ \ol{T} })$ $ \xch $
\xEH
$\xci$ $(\xbm dp)$ (16.2)
\xEH
$f(Y) \xcs (X-f(X)) \xEd \xCQ $ $ \xch $
\xEH
\xEH
\xEP

\cline{3-3}

\xEH
$ \xCN Con( \ol{ \ol{T \xco T' } } \xcv T' )$
\xEH
$\xcI$ without $(\xbm dp)$ (16.3)
\xEH
$f(X \xcv Y) \xcs Y= \xCQ$
\xEH
\xEH
\xEP

\hline

\xEH
$(Log \xcv ' )$
\xEH
$\xch$ $(\xbm\xcc)+(\xbm=)$ (17.1)
\xEH
$( \xbm \xcv ' )$
\xEH
\xEH
\xEP

\cline{3-3}

\xEH
$Con( \ol{ \ol{T' } } \xcv T),$ $ \xCN Con( \ol{ \ol{T' }
} \xcv \ol{ \ol{T} })$ $ \xch $
\xEH
$\xci$ $(\xbm dp)$ (17.2)
\xEH
$f(Y) \xcs (X-f(X)) \xEd \xCQ $ $ \xch $
\xEH
\xEH
\xEP

\cline{3-3}

\xEH
$ \ol{ \ol{T \xco T' } }= \ol{ \ol{T} }$
\xEH
$\xcI$ without $(\xbm dp)$ (17.3)
\xEH
$f(X \xcv Y)=f(X)$
\xEH
\xEH
\xEP

\hline

\xEH
\xEH
\xEH
$( \xbm \xbe )$
\xEH
\xEH
\xEP

\xEH
\xEH
\xEH
$a \xbe X-f(X)$ $ \xch $
\xEH
\xEH
\xEP

\xEH
\xEH
\xEH
$ \xcE b \xbe X.a \xce f(\{a,b\})$
\xEH
\xEH
\xEP

\hline

\end{tabular}

}

\end{turn}

(1) to (7) of the following proposition (in different notation, as the
more
systematic connections were found only afterwards) was already published
in
 \cite{GS08c}, we give it here in totality to complete the picture.
\index{Proposition Ref-Class-Mu-neu}

\bp

$\hspace{0.01em}$

(+++ Orig. No.:  Proposition Ref-Class-Mu-neu +++)

$\hspace{0.01em}$

{\xssc LABEL: {Proposition Ref-Class-Mu-neu}}
\label{Proposition Ref-Class-Mu-neu}

If $f(X)$ is the smallest $ \xCf A$ s.t. $A \xbe \xdf (X),$ then, given
the property on the
left, the one on the right follows.

Conversely, when we define $ \xdf (X):=\{X':f(X) \xcc X' \xcc X\},$ given
the property on
the right, the one on the left follows. For this direction, we assume
that we can use the full powerset of some base set $U$ - as is the case
for
the model sets of a finite language. This is perhaps not too bold, as
we mainly want to stress here the intuitive connections, without putting
too much weight on definability questions.

We assume $ \xCf (iM)$ to hold.

{\footnotesize

\begin{tabular}{|c|c|c|c|}

\hline

(1.1)
\xEH
$(eM\xdi )$
\xEH
$ \xch $
\xEH
$( \xbm wOR)$
\xEP

\cline{1-1}
\cline{3-3}

(1.2)
\xEH
\xEH
$ \xci $
\xEH
\xEP

\hline

(2.1)
\xEH
$(eM\xdi )+(I_\xbo )$
\xEH
$ \xch $
\xEH
$( \xbm OR)$
\xEP

\cline{1-1}
\cline{3-3}

(2.2)
\xEH
\xEH
$ \xci $
\xEH
\xEP

\hline

(3.1)
\xEH
$(eM\xdi )+(I_\xbo )$
\xEH
$ \xch $
\xEH
$( \xbm PR)$
\xEP

\cline{1-1}
\cline{3-3}

(3.2)
\xEH
\xEH
$ \xci $
\xEH
\xEP

\hline

(4.1)
\xEH
$(I \xcv disj )$
\xEH
$ \xch $
\xEH
$( \xbm disjOR)$
\xEP

\cline{1-1}
\cline{3-3}

(4.2)
\xEH
\xEH
$ \xci $
\xEH
\xEP

\hline

(5.1)
\xEH
$(\xdm^+_\xbo) (4)$
\xEH
$ \xch $
\xEH
$( \xbm CM)$
\xEP

\cline{1-1}
\cline{3-3}

(5.2)
\xEH
\xEH
$ \xci $
\xEH
\xEP

\hline

(6.1)
\xEH
$(\xdm^{++})$
\xEH
$ \xch $
\xEH
$( \xbm RatM)$
\xEP

\cline{1-1}
\cline{3-3}

(6.2)
\xEH
\xEH
$ \xci $
\xEH
\xEP

\hline

(7.1)
\xEH
$(I_\xbo )$
\xEH
$ \xch $
\xEH
$( \xbm AND)$
\xEP

\cline{1-1}
\cline{3-3}

(7.2)
\xEH
\xEH
$ \xci $
\xEH
\xEP

\hline

(8.1)
\xEH
$(eM\xdi )+(I_\xbo )$
\xEH
$ \xch $
\xEH
$( \xbm CUT)$
\xEP

\cline{1-1}
\cline{3-3}

(8.2)
\xEH
\xEH
$ \xcI $
\xEH
\xEP

\hline

(9.1)
\xEH
$(eM\xdi )+(I_\xbo )+(\xdm^+_\xbo) (4)$
\xEH
$ \xch $
\xEH
$( \xbm CUM)$
\xEP

\cline{1-1}
\cline{3-3}

(9.2)
\xEH
\xEH
$ \xcI $
\xEH
\xEP

\hline

(10.1)
\xEH
$(eM\xdi )+(I_\xbo )+(eM\xdf )$
\xEH
$ \xch $
\xEH
$( \xbm \xcc \xcd)$
\xEP

\cline{1-1}
\cline{3-3}

(10.2)
\xEH
\xEH
$ \xcI $
\xEH
\xEP

\hline

\end{tabular}

}

\ep

Note that there is no $( \xbm wCM),$ as the conditions $( \xbm  \Xl.)$
imply that the
filter is principal,
and thus that $(I_{ \xbo })$ holds - we cannot ``see'' $ \xCf (wCM)$ alone
with
principal filters.
\index{Proposition Ref-Class-Mu-neu Proof}

\subparagraph{
Proof
}

$\hspace{0.01em}$

(+++ Orig.:  Proof +++)

(1.1) $(eM \xdi )$ $ \xch $ $( \xbm wOR):$

$X-f(X)$ is small in $X,$ so it is small in $X \xcv Y$ by $(eM \xdi ),$ so
$A:=X \xcv Y-(X-f(X)) \xbe \xdf (X \xcv Y),$ but $A \xcc f(X) \xcv Y,$ and
$f(X \xcv Y)$ is the smallest element
of $ \xdf (X \xcv Y),$ so $f(X \xcv Y) \xcc A \xcc f(X) \xcv Y.$

(1.2) $( \xbm wOR)$ $ \xch $ $(eM \xdi ):$

Let $X \xcc Y,$ $X':=Y-$X. Let $A \xbe \xdi (X),$ so $X-A \xbe \xdf (X),$
so $f(X) \xcc X-$A, so
$f(X \xcv X' ) \xcc f(X) \xcv X' \xcc (X-A) \xcv X' $ by prerequisite, so
$(X \xcv X' )-((X-A) \xcv X' )=A \xbe \xdi (X \xcv X' ).$

(2.1) $(eM \xdi )+(I_{ \xbo })$ $ \xch $ $( \xbm OR):$

$X-f(X)$ is small in $X,$ $Y-f(Y)$ is small in $Y,$ so both are small in
$X \xcv Y$ by
$(eM \xdi ),$ so $A:=(X-f(X)) \xcv (Y-f(Y))$ is small in $X \xcv Y$ by
$(I_{ \xbo }),$ but
$X \xcv Y-(f(X) \xcv f(Y)) \xcc A,$ so $f(X) \xcv f(Y) \xbe \xdf (X \xcv
Y),$ so, as $f(X \xcv Y)$ is the smallest
element of $ \xdf (X \xcv Y),$ $f(X \xcv Y) \xcc f(X) \xcv f(Y).$

(2.2) $( \xbm OR)$ $ \xch $ $(eM \xdi )+(I_{ \xbo }):$

Let again $X \xcc Y,$ $X':=Y-$X. Let $A \xbe \xdi (X),$ so $X-A \xbe \xdf
(X),$ so $f(X) \xcc X-$A. $f(X' ) \xcc X',$
so $f(X \xcv X' ) \xcc f(X) \xcv f(X' ) \xcc (X-A) \xcv X' $ by
prerequisite, so
$(X \xcv X' )-((X-A) \xcv X' )=A \xbe \xdi (X \xcv X' ).$

$(I_{ \xbo })$ holds by definition.

(3.1) $(eM \xdi )+(I_{ \xbo })$ $ \xch $ $( \xbm PR):$

Let $X \xcc Y.$ $Y-f(Y)$ is the largest element of $ \xdi (Y),$ $X-f(X)
\xbe \xdi (X) \xcc \xdi (Y)$ by
$(eM \xdi ),$ so $(X-f(X)) \xcv (Y-f(Y)) \xbe \xdi (Y)$ by $(I_{ \xbo }),$
so by ``largest'' $X-f(X) \xcc Y-f(Y),$
so $f(Y) \xcs X \xcc f(X).$

(3.2) $( \xbm PR)$ $ \xch $ $(eM \xdi )+(I_{ \xbo })$

Let again $X \xcc Y,$ $X':=Y-$X. Let $A \xbe \xdi (X),$ so $X-A \xbe \xdf
(X),$ so $f(X) \xcc X-$A, so
by prerequisite $f(Y) \xcs X \xcc X-$A, so $f(Y) \xcc X' \xcv (X-$A), so
$(X \xcv X' )-(X' \xcv (X-A))=A \xbe \xdi (Y).$

Again, $(I_{ \xbo })$ holds by definition.

(4.1) $(I \xcv disj)$ $ \xch $ $( \xbm disjOR):$

If $X \xcs Y= \xCQ,$ then (1) $A \xbe \xdi (X),B \xbe \xdi (Y) \xch A
\xcv B \xbe \xdi (X \xcv Y)$ and
(2) $A \xbe \xdf (X),B \xbe \xdf (Y) \xch A \xcv B \xbe \xdf (X \xcv Y)$
are equivalent. (By $X \xcs Y= \xCQ,$
$(X-A) \xcv (Y-B)=(X \xcv Y)-(A \xcv B).)$
So $f(X) \xbe \xdf (X),$ $f(Y) \xbe \xdf (Y)$ $ \xch $ (by prerequisite)
$f(X) \xcv f(Y) \xbe \xdf (X \xcv Y).$ $f(X \xcv Y)$
is the smallest element of $ \xdf (X \xcv Y),$ so $f(X \xcv Y) \xcc f(X)
\xcv f(Y).$

(4.2) $( \xbm disjOR)$ $ \xch $ $(I \xcv disj):$

Let $X \xcc Y,$ $X':=Y-$X. Let $A \xbe \xdi (X),$ $A' \xbe \xdi (X' ),$
so $X-A \xbe \xdf (X),$ $X' -A' \xbe \xdf (X' ),$
so $f(X) \xcc X-$A, $f(X' ) \xcc X' -A',$ so $f(X \xcv X' ) \xcc f(X)
\xcv f(X' ) \xcc (X-A) \xcv (X' -A' )$ by
prerequisite, so $(X \xcv X' )-((X-A) \xcv (X' -A' ))=A \xcv A' \xbe \xdi
(X \xcv X' ).$

(5.1) $( \xdm^{+}_{ \xbo })$ $ \xch $ $( \xbm CM):$

$f(X) \xcc Y \xcc X$ $ \xch $ $X-Y \xbe \xdi (X),$ $X-f(X) \xbe \xdi (X)$
$ \xch $ (by $( \xdm^{+}_{ \xbo }),$ (4))
$A:=(X-f(X))-(X-Y) \xbe \xdi (Y)$ $ \xch $
$Y-A=f(X)-(X-Y) \xbe \xdf (Y)$ $ \xch $ $f(Y) \xcc f(X)-(X-Y) \xcc f(X).$

(5.2) $( \xbm CM)$ $ \xch $ $( \xdm^{+}_{ \xbo })$

Let $X-A \xbe \xdi (X),$ so $A \xbe \xdf (X),$ let $B \xbe \xdi (X),$ so
$f(X) \xcc X-B \xcc X,$ so by prerequisite
$f(X-B) \xcc f(X).$
As $A \xbe \xdf (X),$ $f(X) \xcc A,$ so $f(X-B) \xcc f(X) \xcc A \xcs
(X-B)=A-$B, and $A-B \xbe \xdf (X-$B), so
$(X-A)-B=X-(A \xcv B)=(X-B)-(A-B) \xbe \xdi (X-$B), so
$( \xdm^{+}_{ \xbo }),$ (4) holds.

(6.1) $( \xdm^{++})$ $ \xch $ $( \xbm RatM):$

Let $X \xcc Y,$ $X \xcs f(Y) \xEd \xCQ.$ If $Y-X \xbe \xdf (Y),$ then
$A:=(Y-X) \xcs f(Y) \xbe \xdf (Y),$ but by
$X \xcs f(Y) \xEd \xCQ $ $A \xcb f(Y),$ contradicting ``smallest'' of
$f(Y).$ So $Y-X \xce \xdf (Y),$ and
by $( \xdm^{++})$ $X-f(Y)=(Y-f(Y))-(Y-X) \xbe \xdi (X),$ so $X \xcs f(Y)
\xbe \xdf (X),$ so $f(X) \xcc f(Y) \xcs X.$

(6.2) $( \xbm RatM)$ $ \xch $ $( \xdm^{++})$

Let $A \xbe \xdf (Y),$ $B \xce \xdf (Y).$ $B \xce \xdf (Y)$ $ \xch $ $Y-B
\xce \xdi (Y)$ $ \xch $ $(Y-B) \xcs f(Y) \xEd \xCQ.$
Set $X:=Y-$B, so $X \xcs f(Y) \xEd \xCQ,$ $X \xcc Y,$ so $f(X) \xcc f(Y)
\xcs X$ by prerequisite.
$f(Y) \xcc A$ $ \xch $ $f(X) \xcc f(Y) \xcs X=f(Y)-B \xcc A-$B.

(7.1) $( \xdi_{ \xbo })$ $ \xch $ $( \xbm AND)$

Trivial.

(7.2) $( \xbm AND)$ $ \xch $ $( \xdi_{ \xbo })$

Trivial.

(8.1) Let $f(X) \xcc Y \xcc X.$ $Y-f(Y) \xbe \xdi (Y) \xcc \xdi (X)$ by
$(eM \xdi ).$ $f(X) \xcc Y$ $ \xch $
$X-Y \xcc X-f(X) \xbe \xdi (X),$ so by $ \xCf (iM)$ $X-Y \xbe \xdi (X).$
Thus by $(I_{ \xbo })$
$X-f(Y)=(X-Y) \xcv (Y-f(Y)) \xbe \xdi (X),$ so $f(Y) \xbe \xdf (X),$ so
$f(X) \xcc f(Y)$ by definition.

(8.2) $( \xbm CUT)$ is too special to allow to deduce $(eM \xdi ).$
Consider $U:=\{a,b,c\},$ $X:=\{a,b\},$ $ \xdf (X)=\{X,\{a\}\},$ $ \xdf
(Z)=\{Z\}$ for all other
$X \xEd Z \xcc U.$ Then $(eM \xdi )$ fails, as $\{b\} \xbe \xdi (X),$ but
$\{b\} \xce \xdi (U).$
$ \xCf (iM)$ and $(eM \xdf )$ hold. We have to check $f(A) \xcc B \xcc A
\xch f(A) \xcc f(B).$ The only
case where it might fail is $A=X,$ $B=\{a\},$ but it holds there, too.

(9.1) By
Fact 14 in  \cite{GS08c}, (6),
we have $( \xbm CM)+( \xbm CUT) \xcj ( \xbm CUM),$ so the result follows
from (5.1) and (8.1).

(9.2) Consider the same example as in (8.2). $f(A) \xcc B \xcc A \xch
f(A)=f(B)$ holds
there, too, by the same argument as above.

(10.1) Let $f(X) \xcc Y,$ $f(Y) \xcc X.$ So $f(X),f(Y) \xcc X \xcs Y,$ and
$X-(X \xcs Y) \xbe \xdi (X),$
$Y-(X \xcs Y) \xbe \xdi (Y)$ by $ \xCf (iM).$ Thus $f(X),f(Y) \xbe \xdf (X
\xcs Y)$ by $(eM \xdf )$ and
$f(X) \xcs f(Y) \xbe \xdf (X \xcs Y)$ by $(I_{ \xbo }).$ So $X \xcs
Y-(f(X) \xcs f(Y)) \xbe \xdi (X \xcs Y),$ so
$X \xcs Y-(f(X) \xcs f(Y)) \xbe \xdi (X), \xdi (Y)$ by $(eM \xdi ),$ so
$(X-(X \xcs Y)) \xcv (X \xcs Y-f(X) \xcs f(Y))$ $=$ $X-f(X) \xcs f(Y) \xbe
\xdi (X)$ by $(I_{ \xbo }),$
so $f(X) \xcs f(Y) \xbe \xdf (X),$ likewise $f(X) \xcs f(Y) \xbe \xdf
(Y),$ so
$f(X) \xcc f(X) \xcs f(Y),$ $f(Y) \xcc f(X) \xcs f(Y),$ and $f(X)=f(Y).$

(10.2) Consider again the same example as in (8.2), we have to show
that $f(A) \xcc B,$ $f(B) \xcc A$ $ \xch $ $f(A)=f(B).$ The only
interesting case is when one
of $A,B$ is $X,$ but not both. Let e.g. $A=X.$ We then have $f(X)=\{a\},$
$f(B)=B \xcc X,$
and $f(X)=\{a\} \xcc B,$ so $B=\{a\},$ and the condition holds.

$ \xcz $
\\[3ex]


\begin{thebibliography}{xxxxxx}

\addcontentsline{toc}{section}{References}


\bibitem[BB94]{BB94}
S.Ben-David, R.Ben-Eliyahu: ``A modal logic for subjective default
reasoning'', Proceedings LICS-94, 1994

\bibitem[GS08c]{GS08c}
D.Gabbay, K.Schlechta, ``Roadmap for preferential logics'',
to appear in: Journal of applied nonclassical logic,
see also hal-00311941, arXiv 0808.3073

\bibitem[Sch90]{Sch90}
K.Schlechta, ``Semantics for Defeasible Inheritance'', in: L.G.Aiello
(ed.), ``Proceedings ECAI 90'', London, 1990, p.594-597

\bibitem[Sch95-1]{Sch95-1}
K.Schlechta: ``Defaults as generalized quantifiers'',
Journal of Logic and Computation, Oxford,
Vol.5, No.4, p.473-494, 1995

\end{thebibliography}
\end{document}